\documentclass[psamsfonts,reqno]{amsart}
\usepackage{amssymb,eucal,graphics,latexsym}

\newcommand{\dnw}{\mathbin{\downarrow}}
\newcommand{\upw}{\mathbin{\uparrow}}

\newcommand{\jirr}{join-ir\-re\-duc\-i\-ble}

\newcommand{\jirry}{join-ir\-re\-duc\-i\-bil\-i\-ty}

\newcommand{\jsd}{join-sem\-i\-dis\-trib\-u\-tive}
\newcommand{\jsdy}{join-sem\-i\-dis\-trib\-u\-tiv\-i\-ty}

\newcommand{\pup}[1]{\textup{(}{#1}\textup{)}}
\newcommand{\rd}[1]{[{#1}]^{\DD}}

\newcommand{\ol}[1]{\overline{#1}}
\newcommand{\Zm}{{\mathbb{Z}}/m{\mathbb{Z}}}
\newcommand{\Zn}{{\mathbb{Z}}/n{\mathbb{Z}}}

\newcommand{\fsi}[1]{\{1,\dots,#1\}}
\newcommand{\fso}[1]{\{0,\dots,#1\}}
\newcommand{\SUB}{\mathbf{SUB}}
\newcommand{\conc}{{{}^{\frown}}}
\newcommand{\es}{\varnothing}

\newcommand{\tr}{\vartriangleleft}

\newcommand{\utr}{\trianglelefteq}
\newcommand{\nutr}{\ntrianglelefteq}
\newcommand{\dtr}{\mathbin{\vartriangleleft\kern-10pt{\lower
3pt\hbox{$\scriptscriptstyle\neq$}}\kern3pt}}
\newcommand{\trd}{\prec}

\newcommand{\set}[1]{\{{#1}\}}
\newcommand{\setm}[2]{\set{{#1}\mid{#2}}}
\newcommand{\seq}[1]{\langle#1\rangle}
\newcommand{\famm}[2]{\seq{{#1}\mid{#2}}}

\newcommand{\St}{\textup{(S)}}
\newcommand{\Ud}{\textup{(U)}}
\newcommand{\Bo}{\textup{(B)}}
\newcommand{\SD}{$\mathrm{(SD_{\vee})}$}
\newcommand{\SDtwo}{$\mathrm{(SD_{\vee}^2)}$}

\newcommand{\Stj}{\textup{(S$_{\mathrm{j}}$)}}
\newcommand{\Udj}{\textup{(U$_{\mathrm{j}}$)}}
\newcommand{\Boj}{\textup{(B$_{\mathrm{j}}$)}}

\newcommand{\fc}{\mathbf{c}}
\newcommand{\fd}{\mathbf{d}}
\newcommand{\fr}{\mathbf{r}}
\newcommand{\fs}{\mathbf{s}}

\newcommand{\hA}{{\widehat{A}}}
\newcommand{\hB}{{\widehat{B}}}
\newcommand{\hC}{{\widehat{C}}}

\newcommand{\into}{\hookrightarrow}

\newcommand{\DD}{\mathbin{D}}
\newcommand{\KK}{\boldsymbol{K}}
\newcommand{\LL}{\boldsymbol{L}}
\newcommand{\eps}{\varepsilon}

\DeclareMathOperator{\J}{J}
\DeclareMathOperator{\Fil}{Fil}

\newcommand{\Co}{\mathbf{Co}}

\newcommand{\bs}{\mathbf{s}}
\newcommand{\bt}{\mathbf{t}}

\newcommand{\vx}{\mathsf{x}}

\numberwithin{equation}{section}

\theoremstyle{plain}
\newtheorem{lemma}{Lemma}[section]
\newtheorem{theorem}[lemma]{Theorem}
\newtheorem{proposition}[lemma]{Proposition}
\newtheorem{corollary}[lemma]{Corollary}
\newtheorem{example}[lemma]{Example}
\newtheorem{claim}{Claim}
\newtheorem*{sclaim}{Claim}

\newtheorem*{stat}{\name}
\newcommand{\name}{testing}

\theoremstyle{definition}
\newtheorem{definition}[lemma]{Definition}
\newtheorem{notation}[lemma]{Notation}
\newtheorem{problem}{Problem}

\theoremstyle{remark}
\newtheorem{remark}[lemma]{Remark}

\newenvironment{all}[1]{\renewcommand{\name}{#1}\begin{stat}}
                        {\end{stat}}

\newcommand{\qedc}{{\qed}~{\rm Claim~{\theclaim}.}}
\newcommand{\qedsc}{{\qed}~{\rm Claim.}}

\newenvironment{cproof}
{\begin{proof}[Proof of Claim.]}
{\qedc\renewcommand{\qed}{}\end{proof}}

\newenvironment{scproof}
{\begin{proof}[Proof of Claim.]}
{\qedsc\renewcommand{\qed}{}\end{proof}}

\begin{document}

\title[Lattices of order-convex sets, I]%
{Sublattices of lattices of order-convex sets, I.\\
The main representation theorem}

\author[M.~Semenova]{Marina Semenova}
\address[M.~Semenova]{Institute of Mathematics of
the Siberian Branch of RAS\\
Acad. Koptyug prosp. 4\\
630090 Novosibirsk\\
Russia}
\email{semenova@math.nsc.ru}

\author[F.~Wehrung]{Friedrich Wehrung}
\address[F.~Wehrung]{LMNO\\
CNRS UMR 6139\\
D\'epartement de Math\'ematiques\\
Universit\'e de Caen\\
14032 Caen Cedex\\
France}
\email{wehrung@math.unicaen.fr}
\urladdr{http://www.math.unicaen.fr/\~{}wehrung}
\date{\today}


\thanks{The first author was partially supported by RFBR grant no.
99-01-00485, by RFBR grant for young scientists no. 01-01-06178, by GA
UK grant no.~162/1999, by GA CR grant no. 201/99, and by INTAS grant no.
YSF: 2001/1-65. The second author was partially supported by the Fund
of Mobility of the Charles University (Prague), by FRVS grant no.
2125, by institutional grant CEZ:J13/98:113200007, and by the
Barrande program.}

\begin{abstract}
For a partially ordered set $P$, we denote by $\Co(P)$ the lattice of
order-convex subsets of $P$. We find three new lattice identities, \St,
\Ud, and \Bo, such that the following result holds.

\begin{all}{Theorem}
Let $L$ be a lattice. Then $L$ embeds into some lattice of the form
$\Co(P)$ if{f} $L$ satisfies \St, \Ud, and \Bo.
\end{all}

Furthermore, if $L$ has an embedding into some $\Co(P)$, then it has
such an embedding that preserves the existing bounds. If
$L$ is finite, then one can take~$P$ finite, with
 \[
 |P|\leq2|\J(L)|^2-5|\J(L)|+4,
 \]
where $\J(L)$ denotes the set of all \jirr\ elements of $L$.

On the other hand, the partially ordered set $P$ can be chosen in such a
way that there are no infinite bounded chains in $P$ and the undirected
graph of the predecessor relation of $P$ is a tree.
\end{abstract}

\maketitle

\section{Introduction}\label{S:Intro}

For a partially ordered set (from now on \emph{poset})
$\seq{P,\utr}$, a subset $X$ of $P$ is \emph{order-convex}, if
$x\utr z\utr y$ and $\set{x,y}\subseteq X$ implies that $z\in X$, for
all $x$, $y$, $z\in P$. The set $\Co(P)$ of all order-convex subsets
of $P$ forms a lattice under inclusion. This lattice is algebraic,
atomistic, and \jsd\ (see Section~\ref{S:BasicC} for the definitions),
thus it is a special example of a \emph{convex geometry}, see
P.\,H.~Edelman~\cite{Ed}, P.\,H.~Edelman and R.~Jamison \cite{EdJa},
or K.\,V.~Adaricheva, V.\,A.~Gorbunov, and V.\,I.~Tumanov~\cite{AGT}.
Furthermore, it is `biatomic' and satisfies the nonexistence of
so-called `zigzags' of odd length on its atoms. Is is proved in G.
Birkhoff and M.\,K. Bennett \cite{BB} that these conditions
\emph{characterize} the lattices of the form
$\Co(P)$.

One of the open problems of \cite{AGT} is the characterization of all
\emph{sublattices} of the lattices of the form $\Co(P)$.

\begin{all}{Problem~3 of \cite{AGT} for $\Co(P)$}
Describe the subclass of those lattices that are embeddable into
finite lattices of the form $\Co(P)$.
\end{all}

In the present paper, we solve completely this problem, not only in the
finite case but also for arbitrary lattices. Our main result
(Theorem~\ref{T:Main}) is that a lattice~$L$ can be embedded into
some lattice of the form $\Co(P)$ if{f} $L$ satisfies three
completely new identities, that we denote by \St, \Ud, and \Bo.
Furthermore, $P$ can be taken either finite in case $L$ is finite, or
\emph{tree-like} (see Theorem~\ref{T:Main2}).

This result is quite surprising, as it yields the unexpected
consequence (see Corollary~\ref{C:Main}) that the class of all lattices
that can be embedded into some $\Co(P)$ is a
\emph{variety}, thus it is closed under homomorphic images. However,
while it is fairly easy (though not completely trivial) to verify
directly that the class is closed under reduced products and
substructures (thus it is a \emph{quasivariety}), we do not know any
direct proof that it is closed under homomorphic images.

One of the difficulties of the
present work is to guess, for a given $L$, which poset~$P$ will solve
the embedding problem for $L$ (i.e., $L$ embeds into $\Co(P)$).
The first natural guess, that consists of using for $P$ the set of
all \jirr\ elements of~$L$, fails, as illustrated by the two examples
of Section~\ref{S:NonAtEx}. We shall construct $P$ \emph{via sequences}
of \jirr\ elements of $L$. In fact, we are able to
embed~$L$ into $\Co(P)$ for two different sorts of posets~$P$:
\begin{itemize}
\item[(1)] $P$ is finite in case $L$ is finite; this is the
construction of Section~\ref{S:Cheshire}.

\item[(2)] $P$ is tree-like (as defined in Section~\ref{S:BasicC});
this is the construction of Section~\ref{S:Koshka}.
\end{itemize}

The two requirements (1) and (2) above can be simultaneously satisfied
in case~$L$ \emph{has no $\DD$-cycle}, see Theorem~\ref{T:Main2}(iii).
However, the finite lattice $\LL$ of Example~\ref{Ex:SINonAt} can
be embedded into some finite $\Co(Q)$, but into no $\Co(R)$, where
$R$ is a finite tree-like poset, see Corollary~\ref{C:CFq-id}. It is
used to produce, in Section~\ref{S:CFq-id}, a
quasi-identity that holds in all $\Co(R)$, where $R$ is finite and
tree-like (or even what we call `crown-free'), but not in all
finite $\Co(P)$.

We conclude the paper by a list of open problems.

\section{Basic concepts}\label{S:BasicC}

A lattice $L$ is \emph{\jsd}, if it satisfies the axiom
 \begin{equation}
 x\vee y=x\vee z\Rightarrow x\vee y=x\vee(y\wedge z),
 \qquad\text{for all }x,\,y,\,z\in L.\tag*{\SD}
 \end{equation}
We denote by $\J(L)$ the set of \jirr\ elements of
$L$. We say that $L$ is \emph{finitely spatial} (resp., \emph{spatial})
if every element of $L$ is a join of \jirr\ (resp., completely
\jirr) elements of $L$.

We say that $L$ is \emph{lower continuous}, if the equality
 \[
 a\vee\bigwedge X=\bigwedge(a\vee X)
 \]
holds, for all $a\in L$ and all downward directed $X\subseteq L$ such
that $\bigwedge X$ exists (where $a\vee X=\setm{a\vee x}{x\in X}$).
It is well known that every dually algebraic
lattice is lower continuous---see Lemma~2.3 in P. Crawley and R.\,P.
Dilworth~\cite{CrDi}, and spatial (thus finitely spatial)---see Theorem
I.4.22 in G. Gierz \emph{et al.}~\cite{Comp} or Lemma~1.3.2 in V.\,A.
Gorbunov~\cite{Gorb}.

For every element $x$ in a lattice $L$, we put
 \[
 \dnw x=\setm{y\in L}{y\leq x};\qquad\upw x=\setm{y\in L}{y\geq x}.
 \]
If $a$, $b$, $c\in L$ such that $a\leq b\vee c$, we
say that the (formal) inequality $a\leq b\vee c$ is a
\emph{nontrivial join-cover}, if $a\nleq b,c$. We say that it is
\emph{minimal in $b$}, if $a\nleq x\vee c$, for all $x<b$, and we say
that it is a \emph{minimal nontrivial join-cover}, if it is a
nontrivial join-cover and it is minimal in both $b$ and $c$.

The \emph{join-dependency} relation $\DD=\DD_L$ (see R. Freese, J.
Je\v{z}ek, and J.\,B. Nation~\cite{FJN}) is defined on the
\jirr\ elements of $L$ by putting
 \[
 p\DD q,\text{ if }p\neq q\text{ and }\exists x\text{ such that }
 p\leq q\vee x\text{ holds and is minimal in }q.
 \]
It is important to observe that $p\DD q$ implies that $p\nleq q$, for
all $p$, $q\in\J(L)$.

For a poset $P$ endowed with a partial ordering $\utr$,
we shall denote by $\tr$ the corresponding strict
ordering. The set of all order-convex subsets of $P$
forms a lattice under inclusion, that we shall denote by $\Co(P)$. The
meet in $\Co(P)$ is the intersection, while the join is given by
 \[
 X\vee Y=X\cup Y\cup\bigcup
 \setm{z\in P}{\exists\seq{x,y}\in(X\times Y)\cup(Y\times X)
 \text{ such that }x\tr z\tr y},
 \]
for all $X$, $Y\in\Co(P)$.
Let us denote by $\prec$ the predecessor relation of $P$. We say
that a \emph{path} of $P$ is a finite sequence
$\fd=\seq{x_0,\dots,x_{n-1}}$ of \emph{distinct} elements of $P$
such that either $x_i\prec x_{i+1}$ or $x_{i+1}\prec x_i$, for all
$i$ with $0\leq i\leq n-2$; if $n>0$, we say that $\fd$ is a path
from $x_0$ to $x_{n-1}$. We say that the path $\fd$ is
\emph{oriented}, if $x_i\prec x_{i+1}$, for all $i$ with
$0\leq i\leq n-2$. We say that $P$ is \emph{tree-like},
if the following properties hold:

\begin{enumerate}
\item for all $a\utr b$ in $P$, there are $n<\omega$ and $x_0$, \dots,
$x_n\in P$ such that $a=x_0\prec x_1\prec\cdots\prec x_n=b$;

\item for all $a$, $b\in P$, there exists at most one path from $a$ to
$b$.
\end{enumerate}

\section{Dually $2$-distributive lattices}\label{S:2DD}

For a positive integer $n$, the identity of $n$-distributivity is
introduced in A.\,P. Huhn \cite{Huhn72}. In this paper we shall only
need the dual of $2$-distributivity, which is the following
identity:
 \begin{equation*}
 a\wedge(x\vee y\vee z)=
 (a\wedge(x\vee y))\vee(a\wedge(x\vee z))\vee(a\wedge(y\vee z)).
 \end{equation*}
We omit the easy proof of the following lemma, that expresses how dual
$2$-distributivity can be read on the \jirr\ elements.

\begin{lemma}\label{L:2Distr}
Let  $L$ be a dually $2$-distributive lattice. For all $p\in\J(L)$
and all $a$, $b$, $c\in L$, if $p\leq a\vee b\vee c$, then either
$p\leq a\vee b$ or $p\leq a\vee c$ or $p\leq b\vee c$.
\end{lemma}

We observe that for finitely spatial $L$, the converse of
Lemma~\ref{L:2Distr} holds.

The following lemma will be used repeatedly throughout the paper.

\begin{lemma}\label{L:ExConj}
Let $L$ be a dually $2$-distributive, complete, lower continuous
lattice. Let $p\in\J(L)$ and let $a$, $b\in L$ such that
$p\leq a\vee b$ and $p\nleq a,b$.
Then the following assertions hold:
\begin{enumerate}
\item There are minimal $x\leq a$ and $y\leq b$ such that
$p\leq x\vee y$.

\item Any minimal $x\leq a$ and $y\leq b$ such that $p\leq x\vee y$
are \jirr.
\end{enumerate}
\end{lemma}

\begin{proof}
(i) Let $X\subseteq\dnw a$ and $Y\subseteq\dnw b$ be chains such that
$p\leq x\vee y$, for all $\seq{x,y}\in X\times Y$. It follows from the
lower continuity of $L$ that
$p\leq\left(\bigwedge X\right)\vee\left(\bigwedge Y\right)$. The
conclusion of (i) follows from a simple application of Zorn's Lemma.

(ii) {}From $p\nleq a,b$ it follows that both $x$ and $y$ are nonzero.
Suppose that $x=x_0\vee x_1$ for some $x_0$, $x_1<x$. It follows from
the minimality assumption on $x$ that $p\nleq x_0\vee y$ and
$p\nleq x_1\vee y$, whence, by Lemma~\ref{L:2Distr},
$p\leq x_0\vee x_1$, thus $p\leq x\leq a$, a contradiction. Hence $x$ is
\jirr.
\end{proof}

For $p$, $a$, $b\in\J(L)$, we say that $\seq{a,b}$ is a
\emph{conjugate pair} with respect to $p$, if $p\nleq a,b$ and $a$
and $b$ are minimal such that $p\leq a\vee b$; we say then that $b$ is
a \emph{conjugate} of $a$ with respect to $p$. Observe that the
latter relation is symmetric in $a$ and $b$, and that it implies that
$p\DD a$ and $p\DD b$.

\begin{notation}\label{No:J(L,a)}
For a lattice $L$ and $p\in\J(L)$, we put
 \[
 \rd{p}=\setm{x\in\J(L)}{p\DD x}.
 \]
\end{notation}

\begin{corollary}\label{C:ExConj}
Let $L$ be a dually $2$-distributive, complete, lower continuous
lattice, and let $p\in\J(L)$. Then every $a\in\rd{p}$
has a conjugate with respect to~$p$.
\end{corollary}

\begin{proof}
By the definition of join-dependency, there exists $c\in L$ such that
$p\leq a\vee c$ and $p\nleq x\vee c$, for all $x<a$. By
Lemma~\ref{L:ExConj}, there are $a'\leq a$ and $b\leq c$ minimal such
that $p\leq a'\vee b$, and both $a'$ and $b$ are \jirr. It follows that
$a'=a$, whence $b$ is a conjugate of $a$ with respect to $p$.
\end{proof}

\section{Stirlitz, Udav, and Bond}\label{S:Identities}

\subsection{The Stirlitz identity \St\ and the axiom \Stj}
\label{Su:Stirlitz}
Let \St\ be the following identity:
 \[
 a\wedge(b'\vee c)=
 (a\wedge b')\vee\bigvee_{i<2}
 \Bigl(a\wedge(b_i\vee c)\wedge
 \bigl((b'\wedge(a\vee b_i))\vee c\bigr)\Bigr),
 \]
where we put $b'=b\wedge(b_0\vee b_1)$.

\begin{lemma}\label{L:CoPhasSt}
The Stirlitz identity \St\ holds in $\Co(P)$, for any poset
$\seq{P,\utr}$.
\end{lemma}

\begin{proof}
Let $A$, $B$, $B_0$, $B_1$, $C\in\Co(P)$ and $a\in A\cap(B'\vee C)$,
where we put $B'=B\cap(B_0\vee B_1)$. Denote by $D$ the right hand
side of the Stirlitz identity calculated with these parameters.
If $a\in B'$ then $a\in A\cap B'\subseteq D$. If $a\in C$ then
$a\in A\cap C\subseteq D$.

Suppose that $a\notin B'\cup C$. There exist $b\in B'$ and
$c\in C$ such that, say, $b\tr a\tr c$. Since $b\in B_0\vee B_1$,
there are $i<2$ and $b'\in B_i$ such that $b'\utr b$, hence
$a\in A\cap(B_i\vee C)$. Furthermore, $b\in B'\cap(A\vee B_i)$, thus
$a\in(B'\cap(A\vee B_i))\vee C$, so $a\in D$.
\end{proof}

\begin{lemma}\label{L:St2DD}
The Stirlitz identity \St\ implies dual $2$-distributivity.
\end{lemma}

\begin{proof}
Take $b_0=x$, $b_1=y$, $b=x\vee y$, and $c=z$.
\end{proof}

Let \SDtwo\ be the following identity:
 \begin{equation}
 x\vee(y\wedge z)=
 x\vee(y\wedge(x\vee(z\wedge(x\vee y)))).\tag*{\SDtwo}
 \end{equation}
It is well known that \SDtwo\ implies \jsdy\ (that is, the axiom \SD),
see, for example, P. Jipsen and H. Rose \cite[page~81]{JiRo}.

\begin{lemma}\label{L:St2SDtwo}
The Stirlitz identity \St\ implies \SDtwo.
\end{lemma}

\begin{proof}
Let $L$ be a lattice satisfying \St, let $x$, $y$, $z\in L$. Set
$y_2=y\wedge(x\vee(z\wedge(x\vee y)))$. Set $a=b_1=y$, $b=z$,
$c=b_0=x$, and $b'=b\wedge(b_0\vee b_1)=z\wedge(x\vee y)$.
Then the following inequalities hold:
 \begin{align*}
 y_2&=y\wedge\Bigl(x\vee\bigl(z\wedge(x\vee y)\bigr)\Bigr)\\
 &=a\wedge\Bigl(\bigl(b\wedge(b_0\vee b_1)\bigr)\vee c\Bigr)\\
 &\leq(a\wedge b')\vee\bigvee_{i<2}
 \Bigl(a\wedge(b_i\vee c)\wedge
 \bigl((b'\wedge(a\vee b_i))\vee c\bigr)\Bigr)\\
 &=(y\wedge z)\vee(y\wedge x)\vee
 \Bigl(y\wedge\bigl((z\wedge y)\vee x\bigr)\Bigr)\\
 &=(y\wedge z)\vee\Bigl(y\wedge\bigl(x\vee(y\wedge z)\bigr)\Bigr)\\
 &=y\wedge\bigl(x\vee(y\wedge z)\bigr)\\
 &\leq x\vee(y\wedge z).
 \end{align*}
This implies that $x\vee y_2\leq x\vee(y\wedge z)$. Since the converse
inequality holds in any lattice, the conclusion follows.
\end{proof}

We now introduce a lattice-theoretical axiom, the \emph{\jirr\
interpretation of \St}, that we will denote by \Stj.

\begin{quote}
For all $a$, $b$, $b_0$, $b_1$, $c\in\J(L)$, the inequalities
$a\leq b\vee c$, $b\leq b_0\vee b_1$, and $a\neq b$ imply that either
$a\leq\ol{b}\vee c$ for some $\ol{b}<b$ or $b\leq a\vee b_i$ and
$a\leq b_i\vee c$ for some $i<2$.
\end{quote}

Throughout the paper we shall make repeated use of the item
(i) of the following statement. Item
(ii) provides a convenient algorithm for verifying
whether a finite lattice satisfies \St.

\begin{proposition}\label{P:JirrStir}
Let $L$ be a lattice. Then the following assertions hold:
\begin{enumerate}
\item If $L$ satisfies \St, then $L$ satisfies~\Stj.

\item If $L$ is complete, lower continuous, finitely spatial, dually
$2$-distributive, and satisfies \Stj, then $L$ satisfies~\St.
\end{enumerate}
\end{proposition}

\begin{proof}
(i) Let $a\leq b\vee c$, $b\leq b_0\vee b_1$,
and $a\ne b$ for some $a$, $b$, $b_0$, $b_1$, $c\in\J(L)$.
Then the element $b'$ of the Stirlitz identity is
$b'=b\wedge(b_0\vee b_1)=b$; observe also that $a\wedge(b\vee c)=a$.
Therefore, applying \St\ yields
 \begin{align*}
 a&=a\wedge(b'\vee c)\\
 &=(a\wedge b')\vee\bigvee_{i<2}
 \Bigl(a\wedge(b_i\vee c)\wedge
 \bigl((b'\wedge(a\vee b_i))\vee c\bigr)\Bigr)\\
 &=(a\wedge b)\vee\bigvee_{i<2}\Bigl(a\wedge(b_i\vee c)\wedge
 \bigl((b\wedge(a\vee b_i))\vee c\bigr)\Bigr).
 \end{align*}
Since $a$ is \jirr, either $a\leq b$ or
$a\leq(b_i\vee c)\wedge((b\wedge(a\vee b_i))\vee c)$ for some $i<2$.
If $a\leq b$ then $a\leq a\vee c$ with $a<b$ (because $a\neq b$).
Suppose that $a\nleq b$. Then
$a\leq(b_i\vee c)\wedge((b\wedge(a\vee b_i))\vee c)\leq b_i\vee c$
for some $i<2$. If $b\not\leq a\vee b_i$, then
$a\leq\ol{b}\vee c$ for $\ol{b}=b\wedge(a\vee b_i)<b$.

(ii) Put $b'=b\wedge(b_0\vee b_1)$, and let $d$ denote the right hand
side of the identity \St. Since $d\leq a\wedge(b'\vee c)$, we must prove
the converse inequality only. Let $a_1\in\J(L)$ with
$a_1\leq a\wedge(b'\vee c)$. Then $a_1\leq a$ and $a_1\leq b'\vee c$.
If $a_1\leq b'$, then $a_1\leq a\wedge b'\leq d$. If $a_1\leq c$, then
$a_1\leq a\wedge c\leq d$.

Suppose now that $a_1\nleq b',c$. Then, by using
Lemma~\ref{L:ExConj}, we obtain that there are minimal $b'_1\leq b'$
and $c_1\leq c$ such that the following inequality holds,
 \begin{equation}\label{Eq:droit}
 a_1\leq b'_1\vee c_1
 \end{equation}
and both $b'_1$ and $c_1$ are \jirr. {}From $a_1\nleq b'$
it follows that $a_1\nleq b'_1$. If $b'_1\leq b_i$ for some $i<2$,
then the inequalities $b'_1\leq b'\wedge b_i\leq b'\wedge(a\vee b_i)$
and $a_1\leq b'_1\vee c_1\leq(b'\wedge(a\vee b_i))\vee c$ hold; but in
this case, we also have $a_1\leq a\wedge(b_i\vee c)$, whence
$a_1\leq d$. Suppose that $b'_1\nleq b_0,b_1$. Then, by
Lemma~\ref{L:ExConj}, there are \jirr\ elements $d_i\leq b_i$, $i<2$,
such that the following inequality
 \begin{equation}\label{Eq:tordu}
 b'_1\leq d_0\vee d_1
 \end{equation}
holds. It follows from \eqref{Eq:droit}, \eqref{Eq:tordu},
$a_1\nleq b'_1$, the minimality of $b'_1$ in \eqref{Eq:droit}, and
\Stj\ that there exists $i<2$ such that $b'_1\leq a_1\vee d_i$ and
$a_1\leq d_i\vee c_1$. Then the following inequalities hold:
 \[
 a_1\leq a\wedge(d_i\vee c_1)\wedge(b'_1\vee c_1)
 \leq a\wedge(b_i\vee c)\wedge((b'\wedge(a\vee b_i))\vee c)
 \leq d.
 \]
In every case, $a_1\leq d$. Since $L$ is finitely spatial, it follows
that $a\wedge(b'\vee c)\leq d$.
\end{proof}

\subsection{The Bond identity \Bo\ and the axiom \Boj}\label{Su:Bond}
Let \Bo\ be the following identity:
 \begin{align*}
 x\wedge(a_0\vee a_1)\wedge(b_0\vee b_1)=&\bigvee_{i<2}
 \Bigl(\bigl(x\wedge a_i\wedge(b_0\vee b_1)\bigr)\vee
 \bigl(x\wedge b_i\wedge(a_0\vee a_1)\bigr)\Bigr)\\
 &\vee\bigvee_{i<2}\bigl(x\wedge(a_0\vee a_1)\wedge(b_0\vee b_1)
 \wedge(a_0\vee b_i)\wedge(a_1\vee b_{1-i})\bigr).
 \end{align*}

\begin{lemma}\label{L:CoPhasBo}
The Bond identity \Bo\ holds in $\Co(P)$, for any poset
$\seq{P,\utr}$.
\end{lemma}

\begin{proof}
Let $X$, $A_0$, $A_1$, $B_0$, $B_1$ be elements of $\Co(P)$. Denote
by $C$ the right hand side of the Bond identity formed from these
elements. Let
$x\in X\cap(A_0\vee A_1)\cap(B_0\vee B_1)$, we prove that $x\in C$. The
conclusion is obvious if $x\in A_0\cup A_1\cup B_0\cup B_1$, so suppose
that $x\notin A_0\cup A_1\cup B_0\cup B_1$.
Since $x\in(A_0\vee A_1)\setminus(A_0\cup A_1)$,
there are $a_0\in A_0$ and $a_1\in A_1$ such that, say,
$a_0\tr x\tr a_1$. Since $x\in(B_0\vee B_1)\setminus(B_0\cup B_1)$,
there are $b_0\in B_0$ and $b_1\in B_1$ such that either
$b_0\tr x\tr b_1$ or $b_1\tr x\tr b_0$. In the first case, $x$ belongs
to $X\cap(A_0\vee A_1)\cap(B_0\vee B_1)\cap
(A_0\vee B_1)\cap(A_1\vee B_0)$, thus to $C$. In the second case, $x$
belongs to $X\cap(A_0\vee A_1)\cap(B_0\vee B_1)\cap
(A_0\vee B_0)\cap(A_1\vee B_1)$, thus again to $C$.
\end{proof}

We now introduce a lattice-theoretical axiom, the \emph{\jirr\
interpretation of \Bo}, that we will denote by \Boj.

\begin{quote}
For all $x$, $a_0$, $a_1$, $b_0$, $b_1\in\J(L)$, the inequalities
$x\leq a_0\vee a_1,b_0\vee b_1$ imply that either $x\leq a_i$ or
$x\leq b_i$ for some $i<2$ or $x\leq a_0\vee b_0$, $a_1\vee b_1$ or
$x\leq a_0\vee b_1$, $a_1\vee b_0$.
\end{quote}

Throughout the paper we shall make repeated use of the item
(i) of the following statement. Item
(ii) provides a convenient algorithm for verifying
whether a finite lattice satisfies \Bo.

\begin{proposition}\label{P:JirrBond}
Let $L$ be a lattice. Then the following assertions hold:
\begin{enumerate}
\item If $L$ satisfies \Bo, then $L$ satisfies \Boj.

\item If $L$ is complete, lower continuous, finitely spatial, dually
$2$-distributive, and satisfies \Boj, then $L$ satisfies \Bo.
\end{enumerate}
\end{proposition}

\begin{proof}
Item (i) is easy to prove by using the
\Bo\ identity and the \jirry\ of $x$.

(ii) Let $u$ (resp., $v$) denote the left (resp., right) hand side of
the identity \Bo. It is obvious that $v\leq u$. Since $L$ is finitely
spatial, in order to prove that $u\leq v$ it is sufficient to prove
that for all $p\in\J(L)$ such that $p\leq u$, the inequality $p\leq v$
holds. This is obvious if either $p\leq a_i$ or $p\leq b_i$ for some
$i<2$, so suppose that $p\nleq a_i,b_i$, for all $i<2$.
Then, by Lemma~\ref{L:ExConj}, there exist
$x_0$, $x_1$, $y_0$, $y_1\in\J(L)$ such that
$x_i\leq a_i$ and $y_i\leq b_i$, for all $i<2$,
while $p\leq x_0\vee x_1$, $y_0\vee y_1$. By assumption,
we obtain that one of the following assertions holds:
 \begin{align*}
 p&\leq(x_0\vee y_0)\wedge(x_1\vee y_1)
 \leq(a_0\vee b_0)\wedge(a_1\vee b_1);\\
 p&\leq(x_0\vee y_1)\wedge(x_1\vee y_0)
 \leq(a_0\vee b_1)\wedge(a_1\vee b_0).
 \end{align*}
In any case, $p\leq v$, which completes the proof.
\end{proof}

\subsection{The Udav identity \Ud\ and the axiom \Udj}\label{Su:Udav}
Let \Ud\ be the following identity:
 \begin{multline*}
 x\wedge(x_0\vee x_1)\wedge(x_1\vee x_2)\wedge(x_0\vee x_2)\\
 =(x\wedge x_0\wedge(x_1\vee x_2))\vee
 (x\wedge x_1\wedge(x_0\vee x_2))\vee
 (x\wedge x_2\wedge(x_0\vee x_1)).
 \end{multline*}

\begin{lemma}\label{L:CoPhasUd}
The Udav identity \Ud\ holds in $\Co(P)$, for any poset
$\seq{P,\utr}$.
\end{lemma}

\begin{proof}
Let $X$, $X_0$, $X_1$, $X_2$ be elements of $\Co(P)$. Denote by $U$
(resp., $V$) the left hand side (resp., right hand side) of the Udav
identity formed from these elements. It is clear that $U$ contains
$V$. Conversely, let $x\in U$, we prove that
$x$ belongs to $V$. This is clear if $x\in X_0\cup X_1\cup X_2$, so
suppose that $x\notin X_0\cup X_1\cup X_2$. Since
$x\in(X_0\vee X_1)\setminus(X_0\cup X_1)$,
there are $x_0\in X_0$ and $x_1\in X_1$ such that, say,
$x_0\tr x\tr x_1$. Since $x\in(X_1\vee X_2)\setminus(X_1\cup X_2)$,
there are $x'_1\in X_1$ and $x_2\in X_2$ such that either
$x'_1\tr x\tr x_2$ or $x_2\tr x\tr x'_1$. But since $x\tr x_1\in X_1$
and $x\notin X_1$, the first possibility is ruled out, whence
$x_2\tr x\tr x'_1$. Since $x\in(X_0\vee X_2)\setminus(X_0\cup X_2)$,
there are $x'_0\in X_0$ and $x'_2\in X_2$ such that either
$x'_0\tr x\tr x'_2$ or $x'_2\tr x\tr x'_0$. The first possibility is
ruled out by $x_2\tr x$ and $x\notin X_2$, while the second possibility
is ruled out by $x_0\tr x$ and $x\notin X_0$. In any case, we obtain a
contradiction.
\end{proof}

As we already did for \St\ and \Bo, we now introduce a
lattice-theoretical axiom, the
\emph{\jirr\ interpretation of \Ud}, that we will denote by \Udj.

\begin{quote}
For all $x$, $x_0$, $x_1$, $x_2\in\J(L)$, the inequalities
$x\leq x_0\vee x_1,x_0\vee\nobreak x_2,x_1\vee\nobreak x_2$ imply that
either $x\leq x_0$ or $x\leq x_1$ or $x\leq x_2$.
\end{quote}

Throughout the paper we shall make repeated use of the item
(i) of the following statement. Item
(ii) provides a convenient algorithm for verifying
whether a finite lattice satisfies \Ud.

\begin{proposition}\label{P:JirrUdav}
Let $L$ be a lattice. Then the following assertions hold:
\begin{enumerate}
\item If $L$ satisfies \Ud, then $L$ satisfies \Udj.

\item If $L$ is complete, lower continuous, finitely spatial, dually
$2$-distributive, and satisfies both \Boj\ and \Udj, then $L$ satisfies
both \Bo\ and \Ud.
\end{enumerate}
\end{proposition}

\begin{proof}
Item (i) is easy to prove by using the
\Ud\ identity and the \jirry\ of $x$.

(ii) We have already seen in Proposition~\ref{P:JirrBond} that $L$
satisfies \Bo.

Let $u$ (resp., $v$) be the left hand side (resp.,
right hand side) of the identity \Ud. It is clear that $v\leq u$.
Let $p\in\J(L)$ such that $p\leq u$, we prove that $p\leq v$. This is
obvious if $p\leq x_i$ for some $i<3$, so suppose that $p\nleq x_i$,
for all $i<3$. Then, by using Lemma~\ref{L:ExConj}, we obtain that
there are \jirr\ elements $p_i$, $p'_i\leq x_i$ ($i<3$) of $L$ such
that the following inequalities hold:
 \begin{equation}\label{Eq:JirrIneq}
 p\leq p_0\vee p_1,\ p'_1\vee p_2,\,p'_0\vee p'_2.
 \end{equation}
Since $p\nleq x_1$ and $L$ satisfies \Boj, it
follows from the first two inequalities of~\eqref{Eq:JirrIneq} that
$p\leq p_0\vee p'_1,p_1\vee p_2$. Similarly, from $p\nleq x_2$, the
last two inequalities of~\eqref{Eq:JirrIneq}, and \Boj, we obtain
the inequalities
$p\leq p'_1\vee p'_2,p'_0\vee p_2$, and from the first and the last
inequality of \eqref{Eq:JirrIneq}, together with $p\nleq x_0$ and \Boj,
we obtain the inequalities $p\leq p_0\vee p'_2,p'_0\vee p_1$. In
particular, we have obtained the inequalities
 \[
 p\leq p_0\vee p'_1,p'_1\vee p'_2,p_0\vee p'_2,
 \]
whence, by the assumption \Udj, $p\leq x_i$ for some $i<3$, a
contradiction.
\end{proof}

\section{First steps together of the identities \St, \Ud, and \Bo}
\label{S:FirstSteps}

\subsection{Udav-Bond partitions}\label{Su:UBpart}
The goal of this subsection is to prove the following partition result
of the sets $\rd{p}$ (see Notation~\ref{No:J(L,a)}).

\begin{proposition}\label{P:UdBo}
Let $L$ be a complete, lower continuous, dually $2$-distributive
lattice that satisfies \Ud\ and \Bo. Then for every $p\in\J(L)$, there
are subsets $A$ and~$B$ of $\rd{p}$ that satisfy the following
properties:
\begin{enumerate}
\item $\rd{p}=A\cup B$ and $A\cap B=\es$.

\item For all $x$, $y\in\rd{p}$, $p\leq x\vee y$ if{f}
$\seq{x,y}\in(A\times B)\cup(B\times A)$.
\end{enumerate}

Moreover, the set $\set{A,B}$ is uniquely determined by these
properties.
\end{proposition}

The set $\set{A,B}$ will be called the \emph{Udav-Bond partition}
(of $\rd{p}$) associated with~$p$. We observe that every conjugate
with respect to $p$ of an element of $A$ (resp., $B$) belongs to $B$
(resp., $A$).

\begin{proof}
If $\rd{p}=\es$ the result is obvious, so suppose that
$\rd{p}\neq\es$. By Lemma~\ref{L:ExConj}, there are $a$, $b\in\rd{p}$
minimal such that $p\leq a\vee b$. We define $A$ and $B$ by the
formulas
 \[
 A=\setm{x\in\rd{p}}{p\leq x\vee b},\qquad
 B=\setm{y\in\rd{p}}{p\leq a\vee y}.
 \]
Let $x\in\rd{p}$. By Corollary~\ref{C:ExConj}, $x$ has a conjugate
with respect to $p$, denote it by~$y$. By Lemma~\ref{L:ExConj}(ii), $y$
is \jirr, thus $y\in\rd{p}$. By applying \Boj\
to the inequalities $p\leq a\vee b,x\vee y$, we obtain that either
$p\leq a\vee x$, $b\vee y$ or $p\leq a\vee y$, $b\vee x$, thus either
$p\leq a\vee x$ or $p\leq b\vee x$. If both inequalities hold
simultaneously, then, since $p\leq a\vee b$ and by \Udj, we obtain
that $p$ lies below either $a$ or $b$ or $x$, a contradiction. Hence
we have established (i).

Let $x$, $y\in\rd{p}$, we shall establish in which case the
inequality $p\leq x\vee y$ holds. Suppose first that $x\in A$ and
$y\in B$. By applying \Boj\ to the inequalities
$p\leq b\vee x,a\vee y$, we obtain that either $p\leq x\vee y$ or
$p\leq b\vee y$. In the second case, $y\in B$, but $y\in A$, a
contradiction by item (i); hence $p\leq x\vee y$.

Now suppose that $x$, $y\in A$. If $p\leq x\vee y$, then, by
applying \Udj\ to the inequalities $p\leq x\vee y,b\vee x,b\vee y$,
we obtain that $p$ lies below either $x$ or $y$ or $b$, a
contradiction. Hence $p\nleq x\vee y$. The conclusion is the same
for $\seq{x,y}\in B\times B$. This concludes the proof of item (ii).

Finally, the uniqueness of $\set{A,B}$ follows easily from items (i)
and (ii).
\end{proof}

\subsection{Choosing orientation with Stirlitz}\label{Su:Orient}

In this subsection we shall investigate further the configuration on
which \Stj\ is based. The following lemma
suggests an `orientation' of the \jirr\ elements in such a
configuration. More specifically, we are trying to embed the given
lattice into $\Co(P)$, for some poset $\seq{P,\utr}$. Attempting to
define $P$ as $\J(L)$, this would suggest to order the elements $a$,
$b$, $b_0$, and $b_1$ by $c\tr a\tr b$ and $b_{1-i}\tr b\tr b_i$.
Although the elements of $P$ will be defined \emph{via finite
sequences} of elements of $\J(L)$, rather than just elements of
$\J(L)$, this idea will be crucial in the construction of
Section~\ref{S:Koshka}.

\begin{lemma}\label{L:OneDir}
Let $L$ be a lattice satisfying~\Stj\ and \Udj. Let $a$, $b$, $b_0$,
$b_1$, $c\in\J(L)$ such that $a\neq b$ and satisfying the inequalities
$a\leq b\vee c$ with $b$ minimal such, and $b\leq b_0\vee b_1$ with
$b\nleq b_0,b_1$. Then the following assertions hold:
\begin{enumerate}
\item The inequalities $b\leq a\vee b_i$ and $a\leq b_i\vee c$
together are equivalent to the single inequality $b\leq b_i\vee c$,
for all $i<2$.

\item There is exactly one $i<2$ such that $b\leq b_i\vee c$.
\end{enumerate}
\end{lemma}

\begin{proof}
We first observe that $b\nleq c$ (otherwise $a\leq c$).
If $b\leq a\vee b_i$ and $a\leq b_i\vee c$, then obviously
$b\leq b_i\vee c$. Suppose, conversely, that $b\leq b_0\vee c$. If
$b\leq b_1\vee c$, then, by observing that $b\leq b_0\vee b_1$
and applying \Udj, we obtain that either
$b\leq b_0$ or $b\leq b_1$ or $b\leq c$, a contradiction. Hence
$b\nleq b_1\vee c$, the uniqueness statement of (ii) follows.
Furthermore, by \Stj, there exists $i<2$ such
that $b\leq a\vee b_i$ and $a\leq b_i\vee c$, whence $b\leq b_i\vee c$,
thus $i=0$. Therefore, $b\leq a\vee b_0$ and $a\leq b_0\vee c$.
\end{proof}

Next, for a conjugate pair $\seq{b,b'}$ of elements of $\J(L)$ with
respect to some element~$a$ of $\J(L)$, we define
 \begin{equation}\label{Eq:DefCabb'}
 C[b,b']=\setm{x\in\J(L)}{b\DD x\text{ and }b\leq b'\vee x}.
 \end{equation}

\begin{notation}\label{Not:SUB}
Let $\SUB$ denote the class of all lattices that satisfy the
identities \St, \Ud, and \Bo.
\end{notation}

Hence $\SUB$ is a variety of lattices. It is \emph{finitely based},
that is, it is defined by finitely many equations.

\begin{lemma}\label{L:Cab}
Let $L$ be a complete, lower continuous, finitely spatial lattice
in $\SUB$. Let $a$, $b\in\J(L)$ such that $a\DD b$.
Then the equality $C[b,b_0]=C[b,b_1]$ holds, for all conjugates
$b_0$ and $b_1$ of $b$ with respect to $a$.
\end{lemma}

\begin{proof}
We prove, for example, that $C[b,b_0]$ is
contained in $C[b,b_1]$. Let $x\in C[b,b_0]$ (so
$b\leq b_0\vee x$), and suppose that $x\notin C[b,b_1]$ (so
$b\nleq b_1\vee x$). By Corollary~\ref{C:ExConj}, $x$ has a conjugate,
say, $y$, with respect to $b$. Since both relations
$a\leq b\vee b_1$ and $b\leq x\vee y$ are minimal nontrivial
join-covers, it follows from Lemma~\ref{L:OneDir} that either
$b\leq b_1\vee x$ or $b\leq b_1\vee y$, but the first possibility does
not hold. Hence the following inequalities hold:
 \begin{equation}\label{Eq:mconf}
 b\leq b_0\vee x,b_1\vee y,x\vee y.
 \end{equation}
Furthermore, by the uniqueness statement of Lemma~\ref{L:OneDir},
$b\nleq b_0\vee y$. Thus, by \Boj\ and the first two inequalities in
\eqref{Eq:mconf} (observe that $b\nleq b_0,b_1,x,y$), we obtain that
$b\leq b_0\vee b_1$. Hence $a\leq b\vee b_0\leq b_0\vee b_1$, whence
$a\leq b\vee b_0,b\vee b_1, b_0\vee b_1$, a contradiction by \Udj.
\end{proof}

For all $a$, $b\in\J(L)$ such that $a\DD b$, there exists, by
Corollary~\ref{C:ExConj}, a conjugate~$b'$ of $b$ with respect to
$a$. By Lemma~\ref{L:Cab}, for fixed $a$, the value of $C[b;b']$ does
not depend of~$b'$. This entitles us to \emph{define}
 \begin{equation}\label{Eq:Cab}
 C(a,b)=C[b,b'],\text{ for any conjugate }b'\text{ of }b
 \text{ with respect to }a.
 \end{equation}

\begin{lemma}\label{L:CabUB}
Let $a$, $b\in\J(L)$ such that $a\DD b$. Then the set
$\set{C(a,b),\rd{b}\setminus C(a,b)}$ is the Udav-Bond partition of
$\rd{b}$ associated with $b$.
\end{lemma}

\begin{proof}
It suffices to prove that the assertions (i) and (ii)
of Proposition~\ref{P:UdBo} are satisfied by the set
$\set{C(a,b),\rd{b}\setminus C(a,b)}$. We first observe the following
immediate consequence of Lemma~\ref{L:OneDir}.

\begin{sclaim}
For any $x\in\rd{b}$ and any conjugate $x'$ of $x$,
$x\notin C(a,b)$ if{f} $x'\in C(a,b)$.
\end{sclaim}

{}From now on we fix a conjugate $b'$ of $b$ with respect to $a$.
Let $x$, $y\in\rd{b}$, let $x'$ (resp., $y'$) be a
conjugate of $x$ (resp., $y$) with respect to $b$.

Suppose first that $x\in C(a,b)$ and $y\notin C(a,b)$, we prove
that $b\leq x\vee y$. It follows from the Claim above that
$y'\in C(a,b)$, whence the inequalities $b\leq b'\vee x,b'\vee y'$ hold,
hence, by \Udj, $b\nleq x\vee y'$. But $b\leq x\vee x',y\vee y'$, thus,
since $b\nleq x,x',y,y'$ and by \Boj, the
inequality $b\leq x\vee y$ holds.

Suppose next that $x$, $y\in C(a,b)$. Since
$b\leq b'\vee x,b'\vee y$, the inequality
$b\leq x\vee y$ would yield, by \Udj, a
contradiction; whence $b\nleq x\vee y$.

Suppose, finally, that $x$, $y\notin C(a,b)$. Thus, by the Claim,
$y'\in C(a,b)$, whence, by the above, $b\leq x\vee y',y\vee y'$,
whence, by \Udj, $b\nleq x\vee y$.
\end{proof}

\subsection{Stirlitz tracks}\label{Su:StTracks}
Throughout this subsection, we shall fix a lattice $L$ satisfying
the identities \St, \Ud, and \Bo. By Lemma~\ref{L:St2DD}, $L$ is
dually $2$-distributive as well. Furthermore, it follows from
Propositions~\ref{P:JirrStir}, \ref{P:JirrBond}, and \ref{P:JirrUdav}
that $L$ satisfies \Stj, \Udj, and \Boj.

\begin{definition}\label{D:StTracks}
For a natural number $n$, a \emph{Stirlitz track} of length $n$ is a
pair $\sigma=\seq{\famm{a_i}{0\leq i\leq n},
\famm{a'_i}{1\leq i\leq n}}$, where the elements $a_i$ for $0\leq
i\leq n$ and $a'_i$ for
$1\leq i\leq n$ are \jirr\ and the following conditions are satisfied:
\begin{enumerate}
\item the inequality $a_i\leq a_{i+1}\vee a'_{i+1}$ holds, for all
$i\in\fso{n-1}$, and it is a minimal nontrivial join-cover;

\item the inequality $a_i\leq a'_i\vee a_{i+1}$ holds, for all
$i\in\fsi{n-1}$.
\end{enumerate}
We shall call $a_0$ the \emph{base} of $\sigma$. Observe that
$a_i\DD a_{i+1}$, for all $i\in\fso{n-1}$.
\end{definition}

Observe that if $\sigma$ is a Stirlitz track as above, then, by
Lemma~\ref{L:OneDir}, the following inequalities also hold:
\begin{align}
a_{i+1}&\leq a_i\vee a_{i+2};\label{Eq:MoreSt1}\\
a_i&\leq a'_{i+1}\vee a_{i+2}\label{Eq:MoreSt2},
\end{align}
for all $i\in\fso{n-2}$.

The main property that we will need about Stirlitz tracks is the
following.

\begin{lemma}\label{L:StTracks}
For a positive integer $n$, let
$\sigma=\seq{\famm{a_i}{0\leq i\leq n},\famm{a'_i}{1\leq i\leq n}}$
be a Stirlitz track of length $n$. Then the inequalities
$a_i\leq a_0\vee a_n$ and $a_i\leq a'_1\vee a_n$ hold, for all
$i\in\fso{n}$. Furthermore, $0\leq k<l\leq n$ implies that
$a_k\nleq a_l$; in particular, the elements $a_i$, for
$0\leq i\leq n$, are distinct.
\end{lemma}

\begin{proof}
We argue by induction on $n$. The result is trivial for $n=1$, and it
follows from \eqref{Eq:MoreSt1} and \eqref{Eq:MoreSt2} for $n=2$. Suppose
that the result holds for $n\geq 2$, and let
$\sigma=\seq{\famm{a_i}{0\leq i\leq n+1},
\famm{a'_i}{1\leq i\leq n+1}}$ be a Stirlitz track of length $n+1$. We
observe that $\sigma_*=\seq{\famm{a_i}{0\leq i\leq n},
\famm{a'_i}{1\leq i\leq n}}$ is a Stirlitz track of length~$n$,
whence, by the induction hypothesis, the following inequalities hold:
 \begin{align}
 a_{n-1}&\leq a_0\vee a_n,\label{Eq:0n-1atn}\\
 a_{n-1}&\leq a'_1\vee a_n.\label{Eq:'1n-1atn}
 \end{align}
We first prove that $a_{n-1}\leq a_0\vee a_{n+1}$. Indeed, suppose
that this does not hold. Hence, \emph{a fortiori}
$a_{n-1}\nleq a_0,a_{n+1}$. Hence, by applying \Boj\ to
\eqref{Eq:MoreSt2} (for $i=n-1$) and \eqref{Eq:0n-1atn} and observing
that $a_{n-1}\nleq a_n,a'_n$, we obtain that
$a_{n-1}\leq a_n\vee a_{n+1}$. Therefore,
$a_{n-1}\leq a_n\vee a_{n+1},a_n\vee a'_n,a'_n\vee a_{n+1}$, a
contradiction by \Udj. Hence, indeed, $a_{n-1}\leq a_0\vee a_{n+1}$.
Consequently, by \eqref{Eq:MoreSt1},
$a_n\leq a_{n-1}\vee a_{n+1}\leq a_0\vee a_{n+1}$. Hence, for
$i\in\fso{n}$, it follows from the induction hypothesis (applied to
$\sigma_*$) that $a_i\leq a_0\vee a_n\leq a_0\vee a_{n+1}$.

The proof of the inequalities $a_i\leq a'_1\vee a_{n+1}$,
for $i\in\fso{n}$, is similar, with~$a_0$ replaced by $a'_1$ and
\eqref{Eq:0n-1atn} replaced by \eqref{Eq:'1n-1atn}.

Finally, let $0\leq k<l\leq n$, and suppose that $a_k\leq a_l$. By
applying the previous result to the Stirlitz track
$\seq{\famm{a_{k+i}}{0\leq i\leq l-k},
\famm{a'_{k+i}}{1\leq i\leq l-k}}$, we obtain that $a_{l-1}\leq
a_k\vee a_l=a_l$, a contradiction. Hence
$a_k\nleq a_l$, in particular, $a_k\neq a_l$.
\end{proof}

\begin{lemma}\label{L:Vcase}
For positive integers $m$, $n>0$, let
 \[
 \sigma=\seq{\famm{a_i}{0\leq i\leq m},\famm{a'_i}{1\leq i\leq m}},
 \quad\tau=\seq{\famm{b_j}{0\leq j\leq n},\famm{b'_j}{1\leq j\leq n}}
 \]
be Stirlitz tracks with the same base $p=a_0=b_0$ and $p\leq a_1\vee
b_1$. Then $a_i,b_j\leq a_m\vee b_n$, for all $i\in\fso{m}$ and
$j\in\fso{n}$.
\end{lemma}

\begin{proof}
Suppose first that the inequality $p\leq a_1\vee b'_1$ holds. Then
$p\leq a_1\vee b'_1,b'_1\vee b_1,b_1\vee a_1$, a contradiction by \Udj.
Hence $p\nleq a_1\vee b'_1$, thus, by applying \Boj\ to the
inequalities $p\leq a_1\vee a'_1,b_1\vee b'_1$, we obtain that
$p\leq a'_1\vee b'_1$.

Furthermore, from Lemma~\ref{L:StTracks} it follows that
$a_i\leq p\vee a_m$, for all
$i\in\fso{m}$, and $b_j\leq p\vee b_n$, for all $j\in\fso{n}$, thus it
suffices to prove that $p\leq a_m\vee b_n$. Again, from
Lemma~\ref{L:StTracks} it follows that
$p\leq a'_1\vee a_m,b'_1\vee b_n$. Suppose that $p\nleq a_m\vee b_n$.
Then $p\nleq a'_1,a_m,b'_1,b_n$, thus, by \Boj, $p\leq a'_1\vee b_n$.
Furthermore, we have seen that $p\leq b'_1\vee b_n$ and
$p\leq a'_1\vee b'_1$. Hence, by \Udj, $p$ lies below either $a'_1$ or
$b'_1$ or $b_n$, a contradiction.
\end{proof}

\section{The small poset associated with a lattice in $\SUB$}
\label{S:Cheshire}

Everywhere in this section before Theorem~\ref{T:Main}, we shall fix a
complete, lower continuous, finitely spatial lattice $L$ in $\SUB$.
For every element $p\in\J(L)$, we denote by $\set{A_p,B_p}$ the
Udav-Bond partition of $\rd{p}$ associated with $p$ (see
Subsection~\ref{Su:UBpart}). We let
$+$ and $-$ be distinct symbols, and we put $R=R_0\cup R_-\cup R_+$,
where $R_0$, $R_-$, and $R_+$ are the sets defined as follows.
 \begin{align*}
 R_0&=\setm{\seq{p}}{p\in\J(L)},\\
 R_+&=\setm{\seq{a,b,+}}{a,b\in\J(L),\ a\DD b},\\
 R_-&=\setm{\seq{a,b,-}}{a,b\in\J(L),\ a\DD b}.
 \end{align*}
We define a map $e\colon R\to\J(L)$ by putting $e(\seq{p})=p$, for
all $p\in\J(L)$, while $e(\seq{a,b,+})=e(\seq{a,b,-})=b$, for all
$a$, $b\in\J(L)$ with $a\DD b$.

Let $\prec$ be the binary relation on $R$ that consists of the
following pairs:
 \begin{align}
 \seq{p,a,-}&\prec\seq{p}\prec\seq{p,b,+}&&\text{whenever }
 a\in A_p\text{ and }b\in B_p,\label{Eq:prec1}\\
 \seq{b,c,-}&\prec\seq{a,b,+}\prec\seq{b,d,+}
 &&\text{and}\label{Eq:prec2}\\
 \seq{b,d,-}&\prec\seq{a,b,-}\prec\seq{b,c,+},&&\text{whenever }
 c\in\rd{b}\setminus C(a,b)\text{ and }d\in C(a,b).\label{Eq:prec3}
 \end{align}

\begin{lemma}\label{L:Ch2St}
Let $\eps\in\set{+,-}$, let $n<\omega$, and let $a_0$, \dots, $a_n$,
$b_0$, \dots, $b_n\in\J(L)$ such that $a_i\DD b_i$, for all
$i\in\fso{n}$ and
$\seq{a_0,b_0,\eps}\prec\cdots\prec\seq{a_n,b_n,\eps}$. Then exactly
one of the following cases occurs:
\begin{enumerate}
\item $\eps=+$ and, putting $a_{n+1}=b_n$, the equality $a_{i+1}=b_i$
holds, for all $i\in\fso{n}$, while there are \jirr\ elements $a'_1$,
\dots, $a'_{n+1}$ of $L$ such that
$\seq{\famm{a_i}{0\leq i\leq n+1},\famm{a'_i}{1\leq i\leq n+1}}$ is a
Stirlitz track.

\item $\eps=-$ and, putting $a_{-1}=b_0$, the equality $a_{i-1}=b_i$
holds, for all $i\in\fso{n}$, while there are \jirr\ elements
$a'_{-1}$, \dots, $a'_{n-1}$ of $L$ such that
$\seq{\famm{a_{n-i}}{0\leq i\leq n+1},
\famm{a'_{n-i}}{1\leq i\leq n+1}}$ is a Stirlitz track.
\end{enumerate}
\end{lemma}

\begin{proof}
Suppose that $\eps=+$ (the proof for $\eps=-$ is similar). We argue by
induction on $n$. If $n=0$, then, from the assumption that $a_0\DD b_0$
and by using Corollary~\ref{C:ExConj}, we obtain a conjugate $a'_1$ of
$b_0$ with respect to $a_0$, and $\seq{\seq{a_0,a_1},\seq{a'_1}}$ is
obviously a Stirlitz track.

Suppose that $n>0$. {}From the assumption that
$\seq{a_{n-1},b_{n-1},+}\prec\seq{a_n,b_n,+}$ and the definition of
$\prec$, we obtain that $a_n=b_{n-1}$. Furthermore, from the induction
hypothesis it follows that there exists a Stirlitz track of the form
 \[
 \seq{\famm{a_i}{0\leq i\leq n},\famm{a'_i}{1\leq i\leq n}}.
 \]
Put $a_{n+1}=b_n$, and let $a'_{n+1}$ be a conjugate of $a_{n+1}$ with
respect to $a_n$. Using again the assumption that
$\seq{a_{n-1},b_{n-1},+}\prec\seq{a_n,b_n,+}$, we obtain the inequality
$a_n\leq a'_n\vee a_{n+1}$. Therefore,
$\seq{\famm{a_i}{0\leq i\leq n+1},\famm{a'_i}{1\leq i\leq n+1}}$ is a
Stirlitz track.
\end{proof}

Let $\utr$ denote the reflexive and transitive closure of $\prec$.

\begin{lemma}\label{L:PartOrd}
The relation $\utr$ is a partial ordering on $R$, and $\prec$ is the
predecessor relation of $\utr$.
\end{lemma}

\begin{proof}
We need to prove that for any $n>0$, if $r_0\prec\cdots\prec r_n$ in
$R$, then $r_0\neq r_n$. We have three cases to consider.

\smallskip

\noindent\textbf{Case 1.} $r_0\in R_+$. In this case,
$r_i=\seq{a_i,b_i,+}\in R_+$, for all $i\in\fsi{n}$. By
Lemma~\ref{L:Ch2St}, if we put $a_{n+1}=b_n$, then $a_{i+1}=b_i$, for
all $i\in\fso{n}$, and there are \jirr\ elements $a'_1$, \dots,
$a'_{n+1}$ of $L$ such that
 \[
 \seq{\famm{a_i}{0\leq i\leq n+1},\famm{a'_i}{1\leq i\leq n+1}}
 \]
is a Stirlitz track. In particular, by Lemma~\ref{L:StTracks},
$a_0\neq a_n$, whence $r_0\neq r_n$.
\smallskip

\noindent\textbf{Case 2.} $r_0\in R_0$. Then $r_i\in R_+$, for all
$i\in\fsi{n}$, thus $r_0\neq r_n$.
\smallskip

\noindent\textbf{Case 3.} $r_0\in R_-$. If $r_n\notin R_-$, then
$r_0\neq r_n$. Suppose that $r_n\in R_-$. Then $r_i=\seq{a_i,b_i,-}$
belongs to $R_-$, for all $i\in\fso{n}$. By Lemma~\ref{L:Ch2St}, if we
put $a_{-1}=b_0$, then $a_{i-1}=b_i$, for all $i\in\fso{n}$, and there
are \jirr\ elements $a'_{-1}$, \dots, $a'_{n-1}$ of $L$ such that
$\seq{\famm{a_{n-i}}{0\leq i\leq n+1},
\famm{a'_{n-i}}{1\leq i\leq n+1}}$ is a Stirlitz track. In particular,
by Lemma~\ref{L:StTracks}, $a_0\neq a_n$, whence $r_0\neq r_n$.
\end{proof}

\begin{definition}\label{D:Isotype}\hfill
\begin{enumerate}
\item Two finite sequences $\fr=\seq{r_0,\dots,r_{n-1}}$ and
$\fs=\seq{s_0,\dots,s_{n-1}}$ of same length of $R$ are
\emph{isotype}, if either $e(r_i)=e(s_i)$, for all $i\in\fso{n-1}$, or
$e(r_i)=e(s_{n-1-i})$, for all $i\in\fso{n-1}$.

\item An oriented path \pup{see Section~\textup{\ref{S:BasicC}}}
$\fr=\seq{r_0,\dots,r_{n-1}}$ of elements of $R$ is
\begin{itemize}
\item[---] \emph{positive} (resp., \emph{negative}), if there are
elements $a_i$, $b_i$ (for $0\leq i<n$) of $\J(L)$ such that
$r_i=\seq{a_i,b_i,+}$ (resp., $r_i=\seq{a_i,b_i,-}$), for all
$i\in\fso{n-1}$.

\item[---] \emph{reduced}, if either it is positive or is has the
form
 \[
 \seq{u_0,\dots,u_{k-1},\seq{p},v_0,\dots,v_{l-1}},
 \]
where $p\in\J(L)$, $\seq{u_0,\dots,u_{k-1}}$ is negative, and
$\seq{v_0,\dots,v_{l-1}}$ is positive.
\end{itemize}

\end{enumerate}

\end{definition}

\begin{lemma}\label{L:ReducedP}
Every oriented path of $R$ is isotype to a reduced oriented path.
\end{lemma}

\begin{proof}
Let $\fr$ be an oriented path of $R$, we prove that $\fr$ is isotype
to a reduced oriented path. If $\fr$ is either positive or reduced
there is nothing to do. Suppose that $\fr$ is neither positive nor
reduced. Then $\fr$ has the form
 \[
 \seq{\seq{a_{k-1},a_k,-},\dots,\seq{a_0,a_1,-},\seq{b_0,b_1,+},
 \dots,\seq{b_{l-1},b_l,+}}
 \]
for some integers $k>0$ and $l\geq 0$. If $l=0$, then $\fr$ is
isotype to the positive path
 \[
 \seq{\seq{a_0,a_1,+},\dots,\seq{a_{k-1},a_k,+}}.
 \]
Suppose now that $l>0$. Since
$\seq{a_0,a_1,-}\prec\seq{b_0,b_1,+}$, two cases can occur.

\begin{itemize}
\item[\textbf{Case 1.}] $a_0=b_1$ and $a_1\notin C(b_0,b_1)$
(see \eqref{Eq:prec2}).
Observe that $\seq{a_0,a_1,-}\prec\seq{a_0}$ if $a_1\in A_{a_0}$
while $\seq{a_0}\prec\seq{a_0,a_1,+}$ if $a_1\in B_{a_0}$
(see \eqref{Eq:prec1}). In the
first case, it follows from Lemma~\ref{L:CabUB} (applied to
$C(a_0,a_1)$) that the sequence
 \[
 \seq{\seq{a_{k-1},a_k,-},\dots,\seq{a_0,a_1,-},\seq{a_0},
 \seq{b_1,b_2,+},\dots,\seq{b_{l-1},b_l,+}}
 \]
is an oriented path, isotype to $\fr$. Similarly, in the second case,
the oriented path
 \[
 \seq{\seq{b_{l-1},b_l,-},\dots,\seq{b_1,b_2,-},\seq{a_0},
 \seq{a_0,a_1,+},\dots,\seq{a_{k-1},a_k,+}}
 \]
is isotype to $\fr$.

\item[\textbf{Case 2.}] $a_1=b_0$ and $b_1\notin C(a_0,a_1)$
(see \eqref{Eq:prec3}).
Observe that $\seq{b_0}\prec\seq{b_0,b_1,+}$ if $b_1\in B_{b_0}$
while $\seq{b_0,b_1,-}\prec\seq{b_0}$ if $b_1\in A_{b_0}$
(see \eqref{Eq:prec1}). In the
first case, the oriented path
 \[
 \seq{\seq{a_{k-1},a_k,-},\dots,\seq{a_1,a_2,-},\seq{b_0},
 \seq{b_0,b_1,+},\dots,\seq{b_{l-1},b_l,+}}
 \]
is isotype to $\fr$. Similarly, in the second case, the oriented path
 \[
 \seq{\seq{b_{l-1},b_l,-},\dots,\seq{b_0,b_1,-},\seq{b_0},
 \seq{a_1,a_2,+},\dots,\seq{a_{k-1},a_k,+}}
 \]
is isotype to $\fr$.
\end{itemize}
This concludes the proof.
\end{proof}

We define a map $\varphi$ from $L$ into the powerset of $R$ as follows:
 \begin{equation}\label{Eq:DefPhi}
 \varphi(x)=\setm{r\in R}{e(r)\leq x},
 \text{ for all }x\in L.
 \end{equation}

\begin{lemma}\label{L:Phi(x)conv}
The set $\varphi(x)$ belongs to $\Co(R,\utr)$, for all $x\in L$.
\end{lemma}

\begin{proof}
It is sufficient to prove that if $r_0\prec\cdots\prec r_n$ in $R$ such
that $e(r_0)$, $e(r_n)\leq x$, the relation $e(r_k)\leq x$
holds whenever $0<k<n$. By Lemma~\ref{L:ReducedP}, it is sufficient
to consider the case where the oriented path
$\fr=\seq{r_0,\dots,r_n}$ is reduced. If it is positive, then, by
Lemma~\ref{L:Ch2St}, there exists a Stirlitz track of the form
 \[
 \seq{\famm{a_i}{0\leq i\leq n+1},\famm{a'_i}{1\leq i\leq n+1}}
 \]
for \jirr\ elements $a_i$, $a'_i$ of $L$ with
$r_i=\seq{a_i,a_{i+1},+}$, for all $i\in\fso{n}$. But then, by
Lemma~\ref{L:StTracks} applied to the Stirlitz track
 \[
 \seq{\famm{a_{i+1}}{0\leq i\leq n},\famm{a'_{i+1}}{1\leq i\leq n}},
 \]
$e(r_k)=a_{k+1}\leq a_1\vee a_{n+1}\leq x$.
Suppose from now on that $\fr$ is not positive. Then three cases can
occur.

\begin{itemize}
\item[\textbf{Case 1.}]
$\fr=\seq{\seq{a_0},\seq{a_0,a_1,+},\dots,\seq{a_{n-1},a_n,+}}$ for
some $a_0$, \dots, $a_n\in\J(L)$. It follows from Lemma~\ref{L:Ch2St}
that there exists a Stirlitz track of the form
 \[
 \seq{\famm{a_i}{0\leq i\leq n},\famm{a'_i}{1\leq i\leq n}},
 \]
hence, by Lemma~\ref{L:StTracks}, $e(r_k)=a_k\leq a_0\vee a_n\leq x$.

\item[\textbf{Case 2.}]
$\fr=\seq{\seq{a_{n-1},a_n,-},\dots,\seq{a_0,a_1,-},\seq{a_0}}$ for
some $a_0$, \dots, $a_n\in\J(L)$. The argument is similar to the one
for Case~1.

\item[\textbf{Case 3.}]
$\fr=\seq{\seq{a_{n'-1},a_{n'},-},\dots,\seq{a_0,a_1,-},\seq{a_0},
\seq{b_0,b_1,+},\dots,\seq{b_{n''-1},b_{n''},+}}$ for some positive
integers $n'$ and $n''$ and \jirr\ $a_0=b_0$, $a_1$, \dots, $a_{n'}$,
$b_1$, \dots, $b_{n''}$.
{}From $\seq{a_0,a_1,-}\prec\seq{a_0}\prec\seq{b_0,b_1,+}$ it follows
that $a_0=b_0\leq a_1\vee b_1$. {}From Lemma~\ref{L:Ch2St} it follows
that there are Stirlitz tracks of the form
 \begin{align*}
 \sigma&=\seq{\famm{a_i}{0\leq i\leq n'},
 \famm{a'_i}{1\leq i\leq n'}},\\
 \tau&=\seq{\famm{b_j}{0\leq j\leq n''},
 \famm{b'_j}{1\leq j\leq n''}},
 \end{align*}
with the same base $a_0=b_0\leq a_1\vee b_1$. Since $e(r_k)$ has
either the form $a_i$, where $0\leq i<n'$, or $b_j$, where
$0\leq j<n''$, it follows from Lemma~\ref{L:Vcase} that
$e(r_k)\leq a_{n'}\vee b_{n''}\leq x$.
\end{itemize}
This concludes the proof.
\end{proof}

\begin{lemma}\label{L:LatEmb}
The map $\varphi$ is a $\seq{0,1}$-lattice embedding from $L$ into
$\Co(R)$.
\end{lemma}

\begin{proof}
It is obvious that $\varphi$ is a $\seq{\wedge,0,1}$-homomorphism.
Let $x$, $y\in L$ such that $x\nleq y$. Since $L$ is finitely spatial,
there exists $p\in\J(L)$ such that $p\leq x$ and $p\nleq y$. Hence,
$\seq{p}\in\varphi(x)\setminus\varphi(y)$, so
$\varphi(x)\nleq\varphi(y)$. Therefore, $\varphi$ is a
$\seq{\wedge,0,1}$-embedding.

Now let $x$, $y\in L$ and let $r\in\varphi(x\vee y)$, we prove
that $r\in\varphi(x)\vee\varphi(y)$. The conclusion is trivial if
$r\in\varphi(x)\cup\varphi(y)$, so suppose
that $r\notin\varphi(x)\cup\varphi(y)$. We need to consider two cases:
\smallskip

\noindent\textbf{Case 1.} $r=\seq{p}$, for some $p\in \J(L)$. So
$p\leq x\vee y$ while $p\nleq x,y$. By Lemma~\ref{L:ExConj}, there are
minimal $u\leq x$ and $v\leq y$ such that $p\leq u\vee v$, hence $u$
and $v$ are \jirr\ and they do not belong to the same side of the
Udav-Bond partition of
$\rd{p}$ associated with $p$ (see Proposition~\ref{P:UdBo}). Hence, by
the definition of $\prec$, either
$\seq{p,u,-}\prec\seq{p}\prec\seq{p,v,+}$ or
$\seq{p,v,-}\prec\seq{p}\prec\seq{p,u,+}$. Since
$\seq{p,u,\eps}\in\varphi(x)$ and $\seq{p,v,\eps}\in\varphi(y)$, for
all $\eps\in\set{+,-}$, it follows from this that
$\seq{p}\in\varphi(x)\vee\varphi(y)$.

\smallskip

\noindent\textbf{Case 2.} $r=\seq{a,b,+}$ for some $a$, $b\in\J(L)$
such that $a\DD b$. So $b\leq x\vee y$ while $b\nleq x,y$. By
Lemma~\ref{L:ExConj}, there are minimal $u\leq x$ and $v\leq y$ such
that $b\leq u\vee v$, hence $u$ and $v$ are \jirr\ and they do not belong
to the same side of the Udav-Bond partition of $\rd{b}$ associated with
$b$ (see Proposition~\ref{P:UdBo}). Hence, it follows from
Lemma~\ref{L:CabUB} that either $u\notin C(a,b)$ and $v\in C(a,b)$ or
$u\in C(a,b)$ and $v\notin C(a,b)$. In the first case,
 \[
 \seq{b,u,-}\prec\seq{a,b,+}\prec\seq{b,v,+},
 \]
while in the second case,
 \[
 \seq{b,v,-}\prec\seq{a,b,+}\prec\seq{b,u,+}.
 \]
Since $\seq{b,u,\eps}\in\varphi(x)$ and
$\seq{b,v,\eps}\in\varphi(y)$, for all $\eps\in\set{+,-}$, it follows
from this that $r\in\varphi(x)\vee\varphi(y)$.

\noindent\textbf{Case 3.} $r=\seq{a,b,-}$ for some $a$, $b\in\J(L)$
such that $a\DD b$.
The proof is similar to the proof of Case~2.
\end{proof}

We can now formulate the main theorem of this paper.

\begin{theorem}\label{T:Main}
Let $L$ be a lattice. Then the following are equivalent:
\begin{enumerate}
\item $L$ embeds into a lattice of the form $\Co(P)$, for some poset
$P$;

\item $L$ satisfies the identities \St, \Ud, and \Bo\ \pup{i.e., it
belongs to the class $\SUB$};

\item $L$ has a lattice embedding into a lattice of the form
$\Co(R)$, for some poset~$R$, that preserves the existing bounds.
Furthermore, if $L$ is finite, then $R$ is finite, with
 \[
 |R|\leq2|\J(L)|^2-5|\J(L)|+4.
 \]
\end{enumerate}
\end{theorem}

\begin{proof}
(i)$\Rightarrow$(ii) follows immediately from Lemmas \ref{L:CoPhasSt},
\ref{L:CoPhasBo}, and \ref{L:CoPhasUd}.

(ii)$\Rightarrow$(iii) Denote by $\Fil L$ the lattice of all \emph{dual
ideals} (= filters) of $L$, ordered by reverse inclusion; if $L$ has no
unit element, then we allow the empty set in $\Fil L$, otherwise we
require filters to be nonempty. This way, $\Fil L$ is complete and the
canonical lattice embedding $x\mapsto\upw x$ from $L$ into $\Fil L$
preserves the existing bounds. It is
well known that $\Fil L$ is a dually algebraic lattice that extends $L$
and that satisfies the same identities as $L$ (see, for example, G.
Gr\"atzer \cite{GLT}), in particular, it belongs to $\SUB$.
Furthermore, $\Fil L$ is dually algebraic, thus lower continuous and
spatial, thus it is \emph{a fortiori} finitely spatial. We consider the
poset $\seq{R,\utr}$ constructed above from $\Fil L$. By Lemmas
\ref{L:Phi(x)conv} and \ref{L:LatEmb}, the canonical map~$\varphi$
defines a $\seq{0,1}$-embedding from $\Fil L$ into $\Co(R)$.

(iii)$\Rightarrow$(i) is trivial.

In case $L$ is finite, put $n=|\J(L)|$, we verify that
$|R|\leq 2n^2-5n+4$ for the poset $\seq{R,\utr}$ constructed above, in
the case where $n\geq2$ (for $n\leq1$ then one can take for $P$ a
singleton). Indeed, it follows from the \jsdy\ of $L$ (that itself
follows from Lemma~\ref{L:St2SDtwo}) that $L$ has at least two
$\DD$-maximal ( = join-prime) elements, hence the number of pairs
$\seq{a,b}$ of elements of $\J(L)$ such that $a\DD b$ is at most
$(n-1)(n-2)$, whence
 \begin{equation}
 |R|\leq2(n-1)(n-2)+n=2n^2-5n+4.\tag*{\qed}
 \end{equation}
\renewcommand{\qed}{}
\end{proof}

\begin{remark}
The upper bound $2|\J(L)|^2-5|\J(L)|+4$ of Theorem~\ref{T:Main}(iii),
obtained for the particular poset $R$ constructed above, is
reached for $L$ defined as the lattice of all order-convex subsets of
a finite chain.
\end{remark}

\begin{corollary}\label{C:Main}
The class of all lattices that can be embedded into some $\Co(P)$
coincides with $\SUB$; it is a finitely based variety. In particular,
it is closed under homomorphic images.
\end{corollary}

Of course, we proved more, for example, the class of all lattices that
can be embedded into some \emph{finite} $\Co(P)$ forms a
\emph{pseudovariety} (see \cite{Gorb}), thus it is closed under
homomorphic images.

\section{The tree-like poset associated with a lattice in $\SUB$}
\label{S:Koshka}

Everywhere in this section before Theorem~\ref{T:Main2}, we shall fix
a complete, lower continuous, finitely spatial lattice $L$ in $\SUB$.
The goal of the present section is to define a tree-like poset
$\Gamma$ and a lattice embedding from $L$ into $\Co(\Gamma)$ that
preserves the existing bounds, see Theorem~\ref{T:Main2}.

The idea to use $\DD$-increasing finite sequences of \jirr\ elements is
introduced in K.\,V. Adaricheva~\cite{Ad1}, where it is proved that
every finite lattice without $\DD$-cycle can be embedded into the
lattice of subsemilattices of some finite meet-semilattice; see also
\cite{AGT}.

We denote by $\Gamma$ the set of all finite, nonempty sequences
$\alpha=\seq{\alpha(0),\ldots,\alpha(n)}$ of elements of $\J(L)$ such
that $\alpha(i)\DD\alpha(i+1)$, for all $i<n$.
We put $|\alpha|=n$ (the \emph{length} of $\alpha$), and we extend
this definition by putting $|\es|=-1$. We further put
$\ol{\alpha}=\seq{\alpha(0),\ldots,\alpha(n-1)}$ (the \emph{truncation}
of $\alpha$) and $e(\alpha)=\alpha(n)$ (the \emph{extremity}
of~$\alpha$). If $\alpha=\ol{\beta}$, we say that $\beta$ is a
\emph{one-step extension} of $\alpha$. Furthermore, for all $n\geq0$,
we put
 \[
 \Gamma_n=\setm{\alpha\in\Gamma}{|\alpha|\leq n},\qquad\text{and}\qquad
 E_n=\Gamma_n\setminus\Gamma_{n-1}\text{ for }n>0.
 \]
For $\alpha\in\Gamma\setminus\Gamma_0$, we say that a \emph{conjugate}
of $\alpha$ is an element $\beta$ of $\Gamma$ such that
$\ol{\alpha}=\ol{\beta}$ and $e(\alpha)$ and $e(\beta)$ are conjugate
with respect to $e(\ol{\alpha})$. It follows from
Corollary~\ref{C:ExConj} that
\emph{every element of $\Gamma\setminus\Gamma_0$ has a conjugate}.
Furthermore, for $\alpha$, $\beta\in\Gamma$, we write
$\alpha\sim\beta$, if either $\alpha=\ol{\beta}$ or
$\beta=\ol{\alpha}$.

For all $n>0$, we define inductively a binary relation $\trd_n$ on
$\Gamma_n$, together with subsets $A_\alpha$ and $B_\alpha$ of
$\rd{e(\alpha)}$ for $\alpha\in\Gamma_{n-1}$.

The induction hypothesis to be satisfied consists of the
following two assertions:
\begin{itemize}
\item[(S1)] $\trd_n$ is acyclic.

\item[(S2)] For all $\alpha$, $\beta\in\Gamma_n$, $\alpha\sim\beta$
if{f} either $\alpha\trd_n\beta$ or $\beta\trd_n\alpha$.
\end{itemize}

For $n=0$, let $\trd_n$ be empty.
\smallskip

The case $n=1$ is the only place where we have some freedom in the
choice of~$\trd_n$. We suppose that we have already used this freedom
for the construction of the poset $\seq{R,\utr}$ of
Section~\ref{S:Cheshire}, that is, for each $p\in\J(L)$, let $A_p$,
$B_p$ such that $\set{A_p,B_p}$ is the Udav-Bond partition of
$\rd{p}$ associated with $p$ (see Subsection~\ref{Su:UBpart}), and we
let $R$ be the poset associated with this choice that we constructed in
Section~\ref{S:Cheshire}. Then we put $A_{\seq{p}}=A_p$ and
$B_{\seq{p}}=B_p$, and we define
 \[
 \trd_1=
 \setm{\seq{\seq{p,a},\seq{p}}}{p\in\J(L),\ a\in A_{\seq{p}}}\cup
 \setm{\seq{\seq{p},\seq{p,b}}}{p\in\J(L),\ b\in B_{\seq{p}}}.
 \]
It is obvious that $\trd_1$ satisfies both (S1) and (S2).
\smallskip

Now suppose having defined $\trd_n$, for $n\geq1$, that satisfies
both (S1) and (S2). For all $\alpha\in E_n$, we define
subsets $A_\alpha$ and $B_\alpha$ of $\rd{e(\alpha)}$ as follows:
\begin{itemize}
\item[\textbf{Case 1.}] $\ol{\alpha}\trd_n\alpha$. Then we put
$A_\alpha=\rd{e(\alpha)}\setminus C(e(\ol{\alpha}),e(\alpha))$ and
$B_\alpha=C(e(\ol{\alpha}),e(\alpha))$.

\item[\textbf{Case 2.}] $\alpha\trd_n\ol{\alpha}$. Then we put
$A_\alpha=C(e(\ol{\alpha}),e(\alpha))$ and
$B_\alpha=\rd{e(\alpha)}\setminus C(e(\ol{\alpha}),e(\alpha))$.
\end{itemize}

Then we define $\trd_{n+1}$ as
 \begin{equation}\label{Eq:Deftrd}
 \trd_{n+1}=\trd_n\cup
 \setm{\seq{\alpha\conc\seq{x},\alpha}}
 {\alpha\in E_n\text{ and }x\in A_\alpha}
 \cup\setm{\seq{\alpha,\alpha\conc\seq{y}}}
 {\alpha\in E_n\text{ and }y\in B_\alpha},
 \end{equation}
where $\seq{\alpha,\beta}\mapsto\alpha\conc\beta$ denotes
\emph{concatenation} of finite sequences.

\begin{lemma}\label{L:ConstrInd}
The relation $\trd_{n+1}$ satisfies both \textup{(S1)} and
\textup{(S2)}.
\end{lemma}

\begin{proof}
It is obvious that $\trd_{n+1}$ satisfies (S2). Now let us prove (S1),
and suppose that $\trd_{n+1}$ has a cycle, say,
$\alpha_0\trd_{n+1}\alpha_1\trd_{n+1}\cdots\trd_{n+1}\alpha_k=\alpha_0$,
where $k\geq2$.
We pick $k$ \emph{minimal} with this property. As
$A_\alpha\cap B_\alpha=\es$, for all $\alpha$, we cannot have $k=2$ as
well, so $k\geq3$.

By the induction hypothesis, one of the elements of the
cycle belongs to $E_{n+1}$, without loss of generality we may assume
that it is the case for $\alpha_0$. Hence, by \eqref{Eq:Deftrd},
$\alpha_1=\ol{\alpha_0}$ belongs to $\Gamma_n$. Let $l$ be the smallest
element of $\fsi{k-1}$ such that $\alpha_{l+1}\notin\Gamma_n$ (it exists
since $\alpha_k=\alpha_0\notin\Gamma_n$). Suppose that $l<k-1$. By (S2)
for~$\trd_{n+1}$, $\alpha_{l+2}=\ol{\alpha_{l+1}}=\alpha_l$, a
contradiction by the minimality of $k$. Hence $l=k-1$, which means that
$\alpha_1$, \ldots, $\alpha_{k-1}\in\Gamma_n$. Hence, since $k-1\geq2$,
we obtain that $\alpha_1\trd_n\cdots\trd_n\alpha_{k-1}=\alpha_1$ is a
$\trd_n$-cycle, a contradiction.
\end{proof}

Lemma~\ref{L:ConstrInd} completes the definition of $\trd_n$, for all
$n>0$. We define $\trd$ as the union over all $n<\omega$ of $\trd_n$.
Hence $\trd$ is an acyclic binary relation on $\Gamma$ such that
$\alpha\sim\beta$ if{f} either $\alpha\trd\beta$ or
$\beta\trd\alpha$, for all $\alpha$, $\beta\in\Gamma$. Since $\trd$
is acyclic, the reflexive and transitive closure $\utr$ of $\trd$ is
a partial ordering on $\Gamma$, for which $\trd$ is exactly the
\emph{predecessor} relation. For the sake of clarity, we rewrite
below the inductive definition of $\trd$ and the sets $A_\alpha$ and
$B_\alpha$ for $\alpha\in\Gamma$.

\begin{itemize}
\item[(a)] For $|\alpha|=0$, $A_\alpha$ and $B_\alpha$ are chosen such
that $\set{A_\alpha,B_\alpha}$ is the Udav-Bond partition of
$\rd{e(\alpha)}$ associated with $e(\alpha)$.

\item[(b)] Suppose that $|\alpha|\geq1$. Then we define $A_\alpha$ and
$B_\alpha$ by
 \[
 \seq{A_\alpha,B_\alpha}=\begin{cases}
 \bigl(\rd{e(\alpha)}\setminus C(e(\ol{\alpha}),e(\alpha)),
 C(e(\ol{\alpha}),e(\alpha))\bigr)&
 \text{ if }\ol{\alpha}\trd\alpha,\\
 \bigl(C(e(\ol{\alpha}),e(\alpha)),
 \rd{e(\alpha)}\setminus C(e(\ol{\alpha}),e(\alpha))\bigr)&
 \text{ if }\alpha\trd\ol{\alpha}.\\
 \end{cases}
 \]
\item[(c)] $\alpha\trd\beta$ implies that $\alpha\sim\beta$.

\item[(d)] $\alpha\conc\seq{x}\trd\alpha$ if{f} $x\in A_\alpha$ and
$\alpha\trd\alpha\conc\seq{x}$ if{f} $x\in B_\alpha$, for all
$\alpha\in\Gamma$ and all $x\in\rd{e(\alpha)}$.
\end{itemize}

By Lemma~\ref{L:CabUB}, the set $\set{A_\alpha,B_\alpha}$ is the
Udav-Bond partition of $\rd{e(\alpha)}$ associated with $\alpha$, for
all $\alpha\in\Gamma$. Therefore, by
Proposition~\ref{P:UdBo} and the definition of
$\trd$, we obtain immediately the following consequence.

\begin{corollary}\label{C:CaUB}
For all $\alpha\in\Gamma$ and all $x$, $y\in\rd{e(\alpha)}$,
$e(\alpha)\leq x\vee y$ if{f} either
$\alpha\conc\seq{x}\trd\alpha\trd\alpha\conc\seq{y}$ or
$\alpha\conc\seq{y}\trd\alpha\trd\alpha\conc\seq{x}$.
\end{corollary}

For $\alpha$, $\beta\in\Gamma$, we denote by $\alpha*\beta$ the largest
common initial segment of $\alpha$ and $\beta$. Observe that
$\alpha*\beta$ belongs to $\Gamma\cup\set{\es}$ and that
$\alpha*\beta=\beta*\alpha$. Put $m=|\alpha|-|\alpha*\beta|$ and
$n=|\beta|-|\alpha*\beta|$. We let $P(\alpha,\beta)$ be the finite
sequence $\seq{\gamma_0,\gamma_1,\ldots,\gamma_{m+n}}$, where the
$\gamma_i$, for $0\leq i\leq m+n$, are defined by $\gamma_0=\alpha$,
$\gamma_{i+1}=\ol{\gamma_i}$, for all $i<m$, $\gamma_{m+n}=\beta$,
and $\gamma_{m+n-j-1}=\ol{\gamma_{m+n-j}}$, for all $j<n$. Hence the
$\gamma_i$-s first decrease from $\gamma_0=\alpha$ to
$\gamma_m=\alpha*\beta$ by successive truncations, then they increase
again from $\gamma_m$ to $\gamma_{m+n}=\beta$ by successive
one-step extensions.
\smallskip

For $\alpha$, $\beta\in\Gamma$, we observe that a path (see
Section~\ref{S:BasicC}) from $\alpha$ to $\beta$ is a finite sequence
$\fc=\seq{\gamma_0,\gamma_1,\ldots,\gamma_k}$
of distinct elements of $\Gamma$ such that
$\gamma_0=\alpha$, $\gamma_k=\beta$, and
$\gamma_i\sim\gamma_{i+1}$, for all $i<k$.

\begin{proposition}\label{P:UniqPath}
For all $\alpha$, $\beta\in\Gamma$, there exists at most one path from
$\alpha$ to $\beta$, and then this path is $P(\alpha,\beta)$.
Furthermore, such a path exists if{f} $\alpha(0)=\beta(0)$.
\end{proposition}

Hence, by using the terminology of Section~\ref{S:BasicC}: the poset
$\seq{\Gamma,\utr}$ is tree-like.

\begin{proof}
Put again $m=|\alpha|-|\alpha*\beta|$ and $n=|\beta|-|\alpha*\beta|$,
and $P(\alpha,\beta)=\seq{\gamma_0,\ldots,\gamma_{m+n}}$. Let
$\fd=\seq{\delta_0,\ldots,\delta_k}$ (for $k<\omega$) be a path from
$\alpha$ to $\beta$. We begin with the following essential observation.

\begin{sclaim}
The path $\fd$ consists of a sequence of truncations followed by a
sequence of one-step extensions.
\end{sclaim}

\begin{scproof}
Suppose that there exists an index $i\in\fsi{k-1}$ such that $\delta_i$
extends both $\delta_{i-1}$ and $\delta_{i+1}$. Then
$\delta_{i-1}=\ol{\delta_i}=\delta_{i+1}$, which contradicts the
fact that all entries of $\fd$ are distinct.

Hence, either $\fd$ consists of a sequence of
truncations, or there exists a least index $l\in\fso{k-1}$ such that
$\delta_{l+1}$ is an extension of $\delta_l$. If $\delta_{i+1}$ is
not an extension of $\delta_i$ for some $i\in\set{l,\ldots,k-1}$,
then, taking the least such $i$, we obtain that $\delta_i$ extends
both $\delta_{i-1}$ and $\delta_{i+1}$, a contradiction by the first
paragraph of the present proof. Hence $\delta_{i+1}$ is a one-step
extension of $\delta_{i}$, for all $i\in\set{l,\ldots, k-1}$.
\end{scproof}

Let $l$ denote the least element of $\fso{k}$ such that $l<k$ implies
that $\delta_{l+1}$ extends $\delta_l$. In particular, $\delta_l$ is a
common initial segment of both $\alpha$ and $\beta$, thus of
$\alpha*\beta$. Furthermore,
 \[
 |\alpha|-l=|\delta_0|-l=|\delta_l|\leq|\alpha*\beta|=|\alpha|-m,
 \]
thus $l\geq m$. Similarly,
 \[
 |\beta|-(k-l)=|\delta_l|\leq|\alpha*\beta|=|\beta|-n,
 \]
thus $k-l\geq n$. In addition, both $\alpha*\beta$ and $\delta_m$ are
initial segments of $\alpha$ of the same length $|\alpha|-m$, thus
$\alpha*\beta=\delta_m$. Similarly, both $\alpha*\beta$ and
$\delta_{k-n}$ are initial segments of $\beta$ of the same length
$|\beta|-n$, whence $\alpha*\beta=\delta_{k-n}$. Therefore,
$\delta_m=\delta_{k-n}$, whence, since all entries of $\fd$ are
distinct, $m=k-n$, so $k=m+n$, whence $l=m$ since $m\leq l\leq k-n$.
It follows then from the claim that $\fd=P(\alpha,\beta)$.

Furthermore, from $\alpha\sim\beta$ it follows that
$\alpha(0)=\beta(0)$, thus the same conclusion follows from the
assumption that there exists a path from $\alpha$ to $\beta$.
Conversely, if $\alpha(0)=\beta(0)$, then $\alpha*\beta$ is nonempty,
thus so are all entries of
$P(\alpha,\beta)$. Hence $P(\alpha,\beta)$ is a path from $\alpha$ to
$\beta$.
\end{proof}

Now we define a map $\pi\colon\Gamma\to R$ by the following rule:
 \[
 \pi(\alpha)=\begin{cases}
 \alpha&\text{if }|\alpha|=0,\\
 \seq{e(\ol{\alpha}),e(\alpha),+}&\text{if }
 \ol{\alpha}\prec\alpha,\\
 \seq{e(\ol{\alpha}),e(\alpha),-}&\text{if }
 \alpha\prec\ol{\alpha},
 \end{cases}
 \qquad\text{for all }\alpha\in\Gamma.
 \]

\begin{lemma}\label{L:piOP}
$\alpha\prec\beta$ in $\Gamma$ implies that
$\pi(\alpha)\prec\pi(\beta)$ in $R$, for all $\alpha$,
$\beta\in\Gamma$. In particular, $\pi$ is order-preserving
\end{lemma}

\begin{proof}
We argue by induction on the
least integer $n$ such that $\alpha$, $\beta\in\Gamma_n$. We need to
consider first the case where $p$, $a$, $b\in\J(L)$, $a\in A_p$,
$b\in B_p$ (so that $\seq{p,a}\prec\seq{p}\prec\seq{p,b}$ in $\Gamma$),
and prove that $\pi(\seq{p,a})\prec\pi(\seq{p})\prec\pi(\seq{p,b})$ in
$R$. But by the definition of $\pi$, the following equalities hold,
 \[
 \pi(\seq{p,a})=\seq{p,a,-},\qquad\pi(\seq{p})=\seq{p},\qquad
 \text{and }\pi(\seq{p,b})=\seq{p,b,+},
 \]
while, by the definition of $\prec$ on $R$,
 \[
 \seq{p,a,-}\prec\seq{p}\prec\seq{p,b,+},
 \]
which solves the case where $n=1$.

The remaining case to consider is where
$\alpha\conc\seq{x}\prec\alpha\prec\alpha\conc\seq{y}$ in $\Gamma$, for
$|\alpha|>0$. Thus $x\in A_\alpha$ and $y\in B_\alpha$, whence
 \begin{align*}
 \pi(\alpha\conc\seq{x})=\seq{e(\alpha),x,-},\\
 \pi(\alpha\conc\seq{y})=\seq{e(\alpha),y,+}.
 \end{align*}
Suppose first that $\ol{\alpha}\prec\alpha$. Then
 \[
 A_\alpha=\rd{e(\alpha)}\setminus C(e(\ol{\alpha}),e(\alpha))
 \text{ while }B_\alpha=C(e(\ol{\alpha}),e(\alpha)).
 \]
Furthermore, $\pi(\alpha)=\seq{e(\ol{\alpha}),e(\alpha),+}$,
while, by the definition of $\prec$ on $R$,
 \[
 \seq{e(\alpha),x,-}\prec\seq{e(\ol{\alpha}),e(\alpha),+}
 \prec\seq{e(\alpha),y,+},
 \]
in other words,
 \[
 \pi(\alpha\conc\seq{x})\prec\pi(\alpha)\prec\pi(\alpha\conc\seq{y}).
 \]
Suppose now that $\alpha\prec\ol{\alpha}$. Then
 \[
 A_\alpha=C(e(\ol{\alpha}),e(\alpha))\text{ while }
 B_\alpha=\rd{e(\alpha)}\setminus C(e(\ol{\alpha}),e(\alpha)).
 \]
Furthermore, $\pi(\alpha)=\seq{e(\ol{\alpha}),e(\alpha),-}$,
while, by the definition of $\prec$ on $R$,
 \[
 \seq{e(\alpha),x,-}\prec\seq{e(\ol{\alpha}),e(\alpha),-}
 \prec\seq{e(\alpha),y,+},
 \]
in other words,
 \[
 \pi(\alpha\conc\seq{x})\prec\pi(\alpha)\prec\pi(\alpha\conc\seq{y}),
 \]
which completes the proof.
\end{proof}

We observe the following immediate consequence of Lemma~\ref{L:piOP}.

\begin{corollary}\label{C:piOP}
One can define a zero-preserving complete meet homomorphism
$\pi^*\colon\Co(R)\to\Co(\Gamma)$ by the rule
 \[
 \pi^*(X)=\pi^{-1}[X],\text{ for all }X\in\Co(R).
 \]
\end{corollary}

We put $\psi=\pi^*\circ\varphi$, where $\varphi\colon L\into\Co(R)$ is
the canonical map defined in Section~\ref{S:Cheshire}. Hence $\psi$ is
a zero-preserving meet homomorphism from $L$ into $\Co(\Gamma)$.
For any $x\in L$, the value $\psi(x)$ is calculated by the same rule
as $\varphi(x)$, see \eqref{Eq:DefPhi}:
 \[
 \psi(x)=\setm{\alpha\in\Gamma}{e(\alpha)\leq x}.
 \]

\begin{lemma}\label{L:gyEmb}
The map $\psi$ is a lattice embedding from $L$ into $\Co(\Gamma)$.
Moreover, $\psi$ preserves the existing bounds.
\end{lemma}

\begin{proof}
The statement about preservation of bounds is obvious. We have
already seen (and it is obvious) that $\psi$ is a meet homomorphism.
Let $x$, $y\in L$ such that $x\nleq y$. Since $L$ is finitely
spatial, there exists $p\in\J(L)$ such that $p\leq x$ and $p\nleq y$;
whence $\seq{p}\in\psi(x)\setminus\psi(y)$. Hence $\psi$ is a meet
embedding from $L$ into $\Co(\Gamma)$.

Let $x$, $y\in L$, let $\alpha\in\psi(x\vee y)$, we prove that
$\alpha\in\psi(x)\vee\psi(y)$. This is obvious if
$\alpha\in\psi(x)\cup\psi(y)$, so suppose that
$\alpha\notin\psi(x)\cup\psi(y)$. Hence $e(\alpha)\leq x\vee y$
while $e(\alpha)\nleq x,y$, thus, by Lemma~\ref{L:ExConj}, there are
minimal $u\leq x$ and $v\leq y$ such that $e(\alpha)\leq u\vee v$, and
both $u$ and $v$ belong to $\rd{e(\alpha)}$. Therefore, by
Corollary~\ref{C:CaUB}, either
$\alpha\conc\seq{u}\trd\alpha\trd\alpha\conc\seq{v}$ or
$\alpha\conc\seq{v}\trd\alpha\trd\alpha\conc\seq{u}$. In both cases,
since $\alpha\conc\seq{u}\in\psi(x)$ and
$\alpha\conc\seq{v}\in\psi(y)$, we obtain that
$\alpha\in\psi(x)\vee\psi(y)$. Therefore, $\psi$ is a join
homomorphism.
\end{proof}

Now we can state the main embedding theorem of the present section.

\begin{theorem}\label{T:Main2}
Let $L$ be a lattice. Then the following assertions are equivalent:
\begin{enumerate}
\item there exists a poset $P$ such that $L$ embeds into $\Co(P)$;

\item $L$ satisfies the identities \St, \Ud, and \Bo\ \pup{i.e., it
belongs to the class $\SUB$};

\item there exists a tree-like
\pup{see Section~\textup{\ref{S:BasicC}}} poset $\Gamma$ such that $L$
has an embedding into $\Co(\Gamma)$ that preserves the existing
bounds. Furthermore, if $L$ is finite without $\DD$-cycle, then
$\Gamma$ is finite.
\end{enumerate}
\end{theorem}

\begin{proof}
(i)$\Rightarrow$(ii) has already been established, see
Theorem~\ref{T:Main}.

(ii)$\Rightarrow$(iii) As in the proof of Theorem~\ref{T:Main}, we
denote by $\Fil L$ the lattice of all filters of $L$, ordered by
reverse inclusion; if $L$ has no unit element, then we allow the empty
set in $\Fil L$, otherwise we require filters to be nonempty. We
consider the poset $\Gamma$ constructed from $\Fil L$ as in
Section~\ref{S:Koshka}. By Lemma~\ref{L:gyEmb},
$L$ embeds into $\Co(\Gamma)$. The finiteness statement of (iii) is
obvious.

(iii)$\Rightarrow$(i) is trivial.
\end{proof}

Even in case $L=\Co(P)$, for a finite totally ordered set $P$, the
poset $\Gamma$ constructed in Theorem~\ref{T:Main2} is not isomorphic
to $P$ as a rule. As it is constructed from finite sequences of
elements of $P$, it does not lend itself to easy graphic
representation. However, many of its properties can be seen on the
simpler poset represented on Figure~5, which is tree-like.

As we shall see in Sections~\ref{S:crowns} and \ref{S:CFq-id}, the
assumption in Theorem~\ref{T:Main2}(iii) that~$L$ be without
$\DD$-cycle cannot be removed.

\section{Non-preservation of atoms}\label{S:NonAtEx}

The posets $R$ and $\Gamma$ that we constructed in Sections
\ref{S:Cheshire} and \ref{S:Koshka} are defined \emph{via} sequences
of \jirr\ elements of $L$. This is to be put in contrast with the main
result of O. Frink \cite{Fr} (see also \cite{GLT}), that embeds any
complemented modular lattice into a geometric lattice: namely, this
construction preserves atoms. Hence the question of the necessity of the
complication of the present paper, that is, using sequences of
\jirr\ elements rather than just \jirr\ elements, is
natural. In the present section we study two examples that show that
this complication is, indeed, necessary.

\begin{example}\label{Ex:LBNonAt}
A finite, atomistic lattice in $\SUB$ without $\DD$-cycle that
cannot be embedded atom-preservingly into any $\Co(T)$.
\end{example}

\begin{proof}
Let $P$ be the nine-element poset represented on the left hand side of
Figure~1, together with order-convex subsets $P_0$, $P_1$, $P_2$,
$Q_0$, $Q_1$, $Q_2$.

\begin{figure}[htb]
\includegraphics{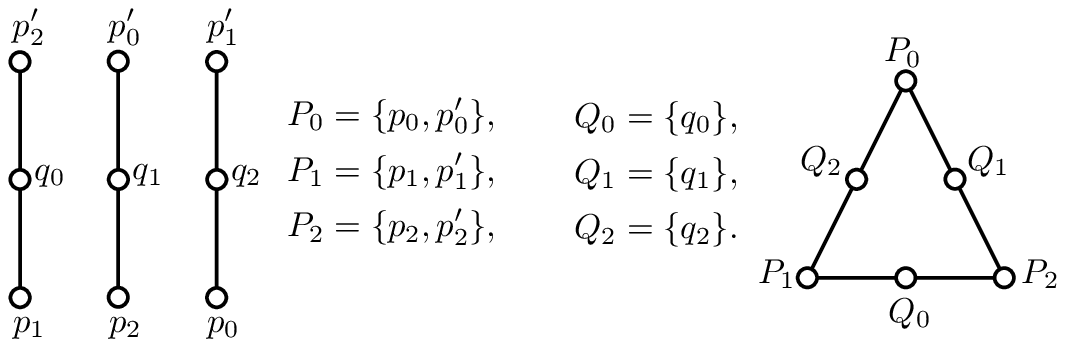}
\caption{The poset $P$ and the geometry of $\KK$.}
\end{figure}

We let $\KK$ be the set of all elements $X$ of $\Co(P)$ such that
$p_i\in X\Leftrightarrow p'_i\in\nobreak X$, for all $i<3$. It is
obvious that $\KK$ is a meet-subsemilattice of $\Co(P)$ which contains
$\set{\es,P}\cup\nobreak\Omega$, where
$\Omega=\set{P_0,P_1,P_2,Q_0,Q_1,Q_2}$. We prove that $\KK$ is a
join-subsemilattice of $\Co(P)$. Indeed, for all $i<3$, both $p_i$
and $p'_i$ are either maximal or minimal in $P$, hence, for all $X$,
$Y\in\Co(P)$, $p_i\in X\vee Y$ if{f} $p_i\in X\cup Y$, and,
similarly, $p'_i\in X\vee Y$ if{f} $p'_i\in X\cup Y$. Hence
$X$, $Y\in\KK$ implies that $X\vee Y\in\KK$.

Therefore, $\KK$ is a sublattice of $\Co(P)$.
It follows immediately that the atoms of $\KK$ are the elements of
$\Omega$, that $\KK$ is atomistic, and the atoms of $\KK$ satisfy the
following relations (see the right half of Figure~1):
 \begin{align*}
 Q_0&\leq P_1\vee P_2;& Q_1&\leq P_0\vee P_2;& Q_2&\leq P_0\vee P_1;\\
 P_0&\nleq P_1\vee P_2;& P_1&\nleq P_0\vee P_2;& P_2&\nleq P_0\vee P_1.
 \end{align*}
Hence, the sequence $P_0P_1P_2P_0P_1$ is a zigzag of length 5 (in the
sense of \cite{BB}). It follows from this and the easy direction of
the main theorem of \cite{BB} that $\KK$ cannot be embedded
atom-preservingly into any $\Co(T)$.
\end{proof}

By contrast, our second example is subdirectly irreducible, but it has
$\DD$-cycles. We shall see in a subsequent paper \cite{SeWe2} that the
latter condition is unavoidable, that is, any finite,
subdirectly irreducible atomistic lattice without $\DD$-cycle that
can be embedded into some $\Co(P)$ can be embedded atom-preservingly
into some finite $\Co(P)$ without $\DD$-cycle.

\begin{example}\label{Ex:SINonAt}
A finite, atomistic, subdirectly irreducible lattice in $\SUB$
that cannot be embedded into $\Co(T)$, for any poset $T$, in an
atom-preserving way.
\end{example}

\begin{proof}
Let $Q$ be the $12$-element poset represented on the left hand side of
Figure~2, together with order-convex subsets $A$, $B$, $C$, $A'$,
$B'$, $C'$.

\begin{figure}[htb]
\includegraphics{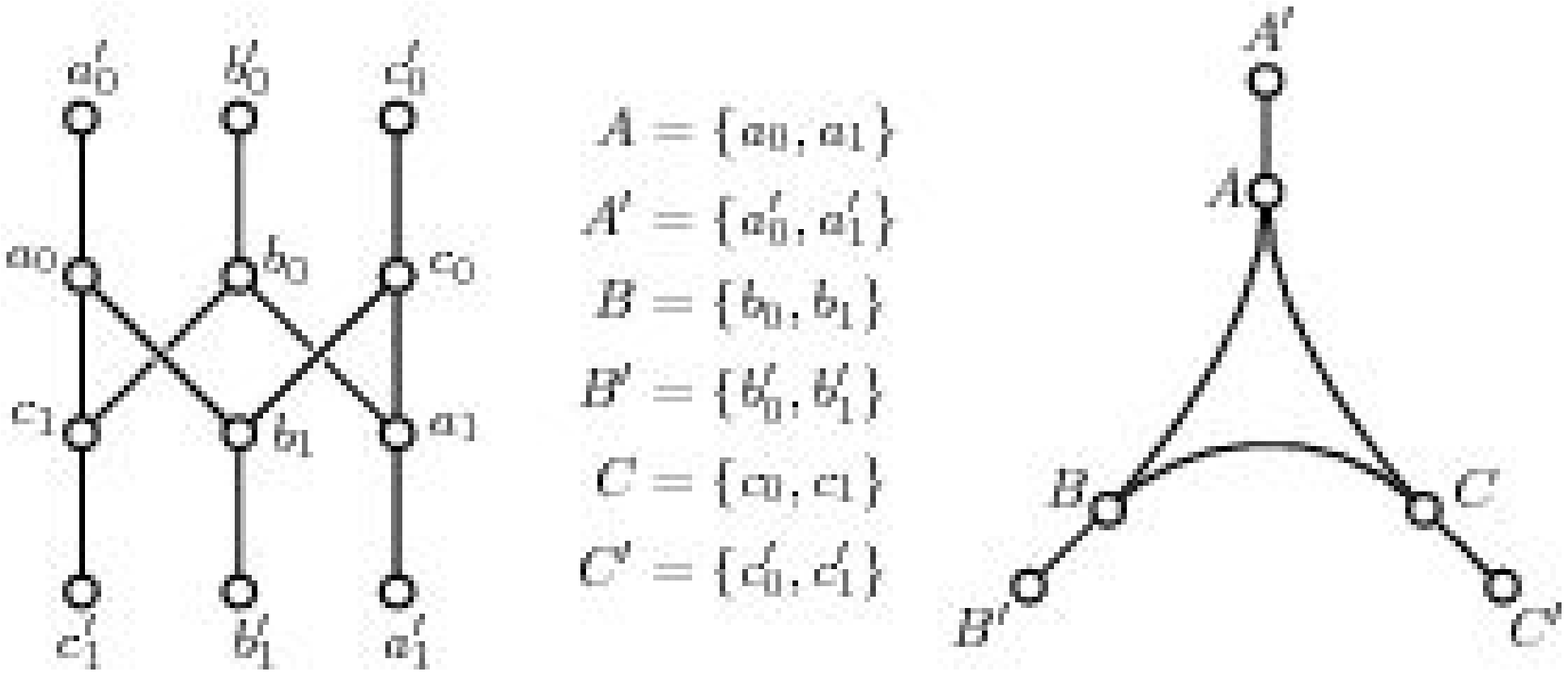}
\caption{The poset $Q$ and the geometry of $\LL$.}
\end{figure}

We let $\sigma$ be the anti-automorphism of $Q$ defined by
$\sigma(a_i)=a_{1-i}$, $\sigma(a'_i)=a'_{1-i}$, $\sigma(b_i)=b_{1-i}$,
$\sigma(b'_i)=b'_{1-i}$, $\sigma(c_i)=c_{1-i}$,
$\sigma(c'_i)=c'_{1-i}$, for all $i<2$, and we let $\LL$ be the set of
all elements $X$ of $\Co(Q)$ such that $\sigma X=X$.
It is obvious that $\LL$ is a meet-subsemilattice of $\Co(Q)$ which
contains $\set{\es,Q}\cup\Omega$, where $\Omega=\set{A,B,C,A',B',C'}$.
We prove that $\LL$ is a join-subsemilattice of $\Co(Q)$. Let $X$,
$Y\in\LL$, we prove that $X\vee Y\in\LL$.

Since both $a'_0$ and $a'_1$ are either maximal or minimal in $Q$, the
equivalence $a'_i\in X\vee Y\Leftrightarrow a'_i\in X\cup Y$ holds,
for all $i<2$, whence $a'_0\in X\vee Y\Leftrightarrow a'_1\in X\vee
Y$. Similarly, $b'_0\in X\vee Y\Leftrightarrow b'_1\in X\vee Y$ and
$c'_0\in X\vee Y\Leftrightarrow c'_1\in X\vee Y$.

Suppose now that $a_0\in X\vee Y$, we prove that $a_1\in X\vee Y$.
If $a_0\in X\cup Y$ this is obvious, so suppose that
$a_0\notin X\cup Y$. Without loss of generality, there are $x\in X$
and $y\in Y$ such that $x\tr a_0\tr y$, whence
$x\in\set{b'_1,b_1,c'_1,c_1}$ and $y=a'_0$. {}From $Y\in\LL$ it
follows that $a'_1\in Y$, thus
$A'\subseteq Y$. Similarly, from $X\in\LL$ it follows that either
$B\subseteq X$ or $C\subseteq X$ or $B'\subseteq X$ or $C'\subseteq X$.
If $B\subseteq X$, then $b_0\in X$, thus, since $a'_1\tr a_1\tr b_0$
and $a'_1\in Y$, we obtain that $a_1\in X\vee Y$. If $B'\subseteq X$,
then $b'_0\in X$, thus, since $a'_1\tr a_1\tr b'_0$ and $a'_1\in Y$, we
obtain again that $a_1\in X\vee Y$. Similar results hold for either
$C\subseteq X$ or $C'\subseteq X$. Therefore, $a_0\in X\vee Y$ implies
that $a_1\in X\vee Y$. By symmetry, we obtain the converse. Similarly,
$b_0\in X\vee Y\Leftrightarrow b_1\in X\vee Y$ and
$c_0\in X\vee Y\Leftrightarrow c_1\in X\vee Y$. Therefore, $X\vee Y$
belongs to $\LL$, which completes the proof that $\LL$ is a sublattice of
$\Co(Q)$.

It follows immediately that the atoms of $\LL$ are the elements of
$\Omega$, that $\LL$ is atomistic, and the atoms of $\LL$ satisfy the
following relations:
 \begin{align*}
 A,B&\leq A'\vee B';& A&\leq A'\vee B;&B\leq A\vee B';\\
 B,C&\leq B'\vee C';& B&\leq B'\vee C;&C\leq B\vee C';\\
 A,C&\leq A'\vee C';& A&\leq A'\vee C;&C\leq A\vee C'.
 \end{align*}
Hence, $\LL$ is subdirectly irreducible, with monolith (i.e., smallest
nonzero congruence) the smallest congruence $\Theta(\es,A)$
identifying $\es$ and $A$, also equal to
$\Theta(\es,B)$ and to $\Theta(\es,C)$. Furthermore, the sequence
$A'B'C'A'B'$ is a zigzag of length 5 (in the sense of \cite{BB}). It
follows from this and the easy direction of the main theorem
of~\cite{BB} that $\LL$ cannot be embedded atom-preservingly into
any~$\Co(T)$.
\end{proof}

\section{Crowns in posets}\label{S:crowns}

We first recall the following classical definition.

\begin{definition}\label{D:crown}
For an integer $n\geq2$, we denote by $\Zn$ the set of integers
modulo $n$. The \emph{$n$-crown} $C_n$ is the poset
with underlying set $(\Zn)\times\set{0,1}$ and ordering defined by
$(i,0),(i+1,0)<(i,1)$, for all $i\in\Zn$.
\end{definition}

The crown $C_n$ is illustrated on Figure~3.

\begin{figure}[htb]
\includegraphics{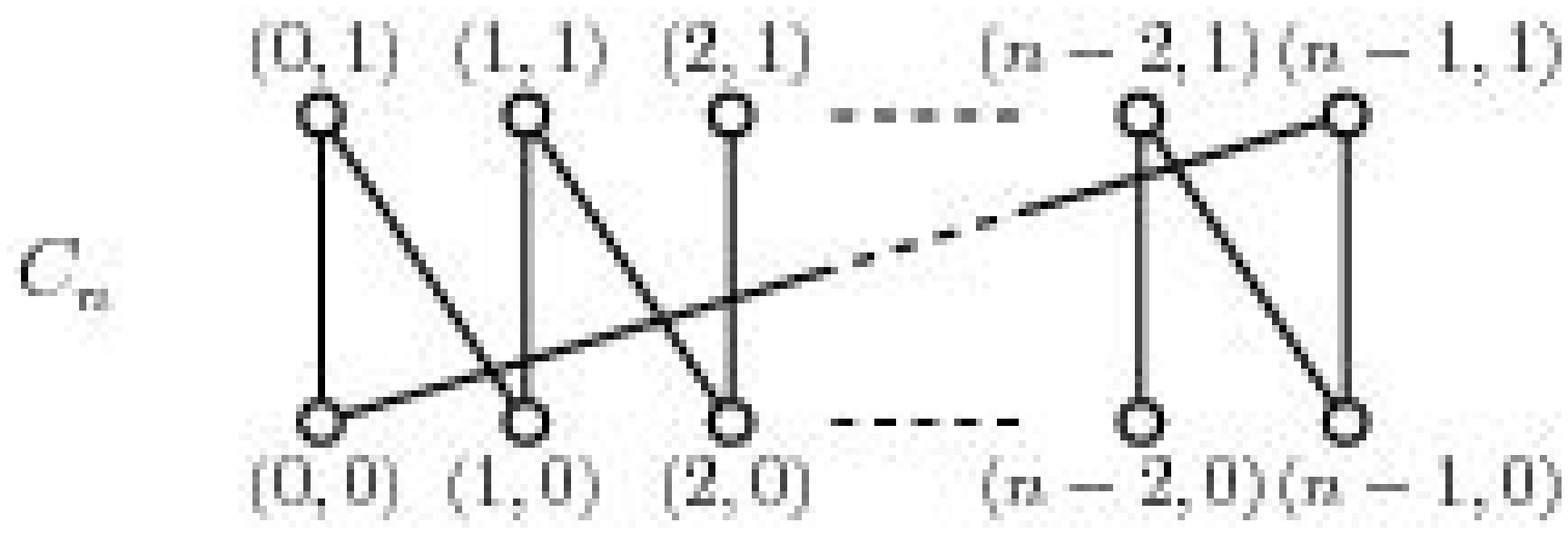}
\caption{The crown $C_n$.}
\end{figure}

We shall mostly deal with sub-crowns of posets.

\begin{definition}\label{D:inncr}
For $n\geq2$ and a poset $(T,\utr)$, a \emph{$n$-crown} of $T$ is a
finite sequence $\famm{\seq{a_i,b_i}}{i\in\Zn}$ of elements of
$T\times T$ such that there exists an order-embedding
$f\colon C_n\into T$ with $f(i,0)=a_i$ and $f(i,1)=b_i$, for all
$i\in\Zn$.
\end{definition}

We shall sometimes identify an integer modulo $n$ with its unique
representative in $\set{0,1,\dots,n-1}$ and a $n$-crown
$\famm{\seq{a_i,b_i}}{i\in\Zn}$ with the finite sequence
 \[
 \seq{\seq{a_0,b_0},\seq{a_1,b_1},\dots,\seq{a_{n-1},b_{n-1}}}.
 \]
The following lemma makes it possible to identify crowns within posets.

\begin{lemma}\label{L:subcrown}
Let $(T,\utr)$ be a poset, let $n\geq3$, and let $a_i$, $b_i$
\pup{$i\in\Zn$} be elements of~$T$. Then the following are equivalent:
\begin{enumerate}
\item $\famm{\seq{a_i,b_i}}{i\in\Zn}$ is a $n$-crown.

\item $a_i\utr b_j$ if{f} $i\in\set{j,j+1}$, for all $i$, $j\in\Zn$.
\end{enumerate}
\end{lemma}

\begin{proof}
(i)$\Rightarrow$(ii) is trivial. Conversely, suppose (ii) satisfied,
we prove that $f\colon C_n\into\nobreak T$ defined by $f(i,0)=a_i$ and
$f(i,1)=b_i$, for all $i\in\Zn$, is an order-embedding. We need to
prove the following assertions:
\begin{enumerate}
\item $a_i\utr a_j$ implies that $i=j$, for all $i$, $j\in\Zn$.
Indeed, if $a_i\utr a_j$, then $a_i\utr b_j,b_{j-1}$ (because
$a_j\utr b_j,b_{j-1}$), thus, by assumption,
$i\in\set{j,j+1}\cap\set{j,j-1}=\set{j}$ (we use here the inequality
$n\geq3$), that is, $i=j$.

\item $b_i\utr b_j$ implies that $i=j$, for all $i$, $j\in\Zn$. The
proof is similar to the one of (i).

\item $b_j\utr a_i$ occurs for no $i$, $j\in\Zn$.
Indeed, suppose that $b_j\utr a_i$. Then $b_j\utr b_i,b_{i-1}$
(because $a_i\utr b_i,b_{i-1}$), thus, by (ii), $j=i=i-1$, a
contradiction.
\end{enumerate}
This concludes the proof.
\end{proof}

\begin{definition}\label{D:CrFree}
A poset $T$ is \emph{crown-free}, if it has no $n$-crown for any
$n\geq3$.
\end{definition}

Strictly speaking, the $2$-crown $C_2$ is crown-free since we are
requiring $n\geq3$ in the definition above. The motivation why we are
putting this slight restriction on~$n$ lies in the following
observation. First, the poset of Figure~4(i) is tree-like, but it
contains the $2$-crown represented on Figure~4(ii); observe also that
the $n$-crown, for any $n\geq2$, is never tree-like.

\begin{figure}[htb]
\includegraphics{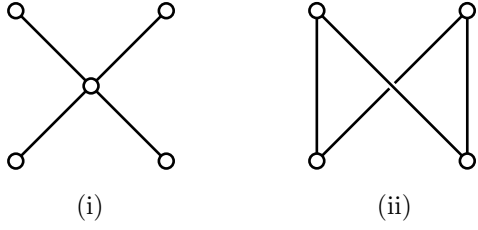}
\caption{A tree-like poset which contains the crown $C_2$.}
\end{figure}

On the other hand, we shall now prove the following result.

\begin{proposition}\label{P:TrCf}
Every tree-like poset is crown-free.
\end{proposition}

As witnessed by the square $\mathbf{2}^2$, the converse of
Proposition~\ref{P:TrCf} does not hold.

\begin{proof}
Let $(T,\utr)$ be a tree-like poset. For $x$, $y\in T$, we denote by
$d(x,y)$ the length of the unique path from $x$ to $y$ if there is
such a path, $\infty$ otherwise. Observe that $x\utr y$ implies that
$d(x,y)<\infty$ (but the converse does not hold as a rule), and then the
unique path from $x$ to $y$ is oriented (see Section~\ref{S:BasicC}).

For a $n$-crown $\gamma=\famm{\seq{a_i,b_i}}{i\in\Zn}$ in $T$, we put
 \[
 \ell(\gamma)=\sum_{i\in\Zn}d(a_i,b_i).
 \]
Suppose that $T$ has a $n$-crown, for some integer $n\geq3$. We pick
such a crown $\gamma=\famm{\seq{a_i,b_i}}{i\in\Zn}$ with
$\ell(\gamma)$ \emph{minimum}. For all $i\in\Zn$, we let
 \begin{align*}
 a_i=x_{i,0}&\prec x_{i,1}\prec\cdots\prec x_{i,p_i}=b_i,\\
 a_{i+1}=y_{i,0}&\prec y_{i,1}\prec\cdots\prec y_{i,q_i}=b_i
 \end{align*}
be the paths from $a_i$ (resp., $a_{i+1}$) to $b_i$, where $\prec$
denotes the predecessor relation of $T$.

\setcounter{claim}{0}
\begin{claim}\label{Cl:noxy1}
$\setm{x_{i,p}}{0\leq p<p_i}\cap\setm{y_{i,q}}{0\leq q<q_i}=\es$,
for all $i\in\Zn$.
\end{claim}

\begin{cproof}
Suppose, to the contrary, that $x_{i,p}=y_{i,q}$ for some
$p\in\fso{p_i-\nobreak 1}$ and $q\in\fso{q_i-1}$. We put $b'_j=b_j$,
for all $j\neq i$ in $\Zn$, while $b'_i=x_{i,p}$. Since
$a_i,a_{i+1}\utr b'_i$, the condition $k\in\set{l,l+1}$ implies that
$a_k\utr b'_l$, for all $k$, $l\in\Zn$. Conversely, let $k$, $l\in\Zn$
such that $a_k\utr b'_l$. {}From $b'_l\utr b_l$ it follows that
$a_k\utr b_l$, whence $k\in\set{l,l+1}$. By Lemma~\ref{L:subcrown},
the family $\gamma'=\famm{\seq{a_k,b'_k}}{k\in\Zn}$ is a $n$-crown.
However,
 \[
 \ell(\gamma')\leq\ell(\gamma)-(p_i-p)<\ell(\gamma),
 \]
which contradicts the minimality of $\ell(\gamma)$.
\end{cproof}

The proof of the following claim is symmetric.

\begin{claim}\label{Cl:noxy2}
$\setm{y_{i,q}}{0<q\leq q_i}\cap\setm{x_{i+1,p}}{0<p\leq p_{i+1}}=\es$,
for all $i\in\Zn$.
\end{claim}

We define a \emph{walk} of $T$ to be a finite sequence
$\fc=\seq{c_0,c_1,\cdots,c_m}$ of elements of~$T$ such that either
$c_i\prec c_{i+1}$ or $c_{i+1}\prec c_i$, for all $i<m$, we say then
that $\fc$ is a walk from $c_0$ to $c_m$. Hence, a nonempty path of
$T$ is a walk with all distinct entries.

Now we let $\fd$ be the finite sequence defined by
 \begin{multline*}
 \fd=\famm{x_{0,k}}{0\leq k\leq p_0}\conc
 \famm{y_{0,q_0-l}}{0<l<q_0}\conc
 \famm{x_{1,k}}{0\leq k\leq p_1}\conc\cdots\\
 \cdots\conc\famm{x_{n-1,k}}{0\leq k\leq p_{n-1}}.
 \end{multline*}
It is obvious that $\fd$ is a walk from $x_{0,0}=a_0$ to
$x_{n-1,p_{n-1}}=b_{n-1}$. We shall now prove that $\fd$ is a path.

Suppose, indeed, that $\fd$ is not a path. Then one of the following
cases occurs.

\begin{itemize}
\item[\textbf{Case 1.}] There are distinct $i$, $j\in\Zn$, together
with $k\in\fso{p_i}$ and $l\in\fso{p_j}$, such that $x_{i,k}=x_{j,l}$.
Then $a_i\utr x_{i,k}=x_{j,l}\utr b_j$, thus $i\in\set{j,j+1}$, while
$a_j\utr x_{j,l}=x_{i,k}\utr b_i$, thus $j\in\set{i,i+1}$. Since
$n\geq3$, we obtain that $i=j$, a contradiction.

\item[\textbf{Case 2.}] There are distinct $i$,
$j\in(\Zn)\setminus\set{n-1}$, together
with $k\in\fsi{q_i-1}$ and $l\in\fsi{q_j-1}$, such that
$y_{i,k}=y_{j,l}$. Then $a_{i+1}\utr y_{i,k}=y_{j,l}\utr b_j$, thus
$i\in\set{j,j-1}$, while $a_{j+1}\utr y_{j,l}=y_{i,k}\utr b_i$, thus
$j\in\set{i,i-1}$, whence, since $n\geq3$, $i=j$, a contradiction.

\item[\textbf{Case 3.}] There are $i\in\Zn$ and
$j\in(\Zn)\setminus\set{n-1}$, together with $k\in\fso{p_i}$ and
$l\in\fsi{q_j-1}$, such that $x_{i,k}=y_{j,l}$. Then from
Claim~\ref{Cl:noxy1} it follows that $i\neq j$, while from
Claim~\ref{Cl:noxy2} it follows that $i\neq j+1$. On the other hand,
$a_i\utr x_{i,k}=y_{j,l}\utr b_j$, thus $i\in\set{j,j+1}$, a
contradiction.
\end{itemize}

Therefore, we have proved that $\fd$ is, indeed, a path from $a_0$ to
$b_{n-1}$. However, the finite sequence
 \[
 \fd'=\famm{y_{n-1,l}}{0\leq l\leq q_{n-1}}
 \]
is a path from $y_{n-1,0}=a_n=a_0$ (the indices are modulo $n$) to
$y_{n-1,q_{n-1}}=b_{n-1}$, thus, by the uniqueness of the path from
$a_0$ to $b_{n-1}$, $\fd=\fd'$. Thus every entry $x$ of
$\fd$ satisfies that $x\utr b_{n-1}$, in particular,
$b_0=x_{0,p_0}\utr b_{n-1}$, a contradiction since $n\neq1$.
\end{proof}

\section{A quasi-identity for $\Co(T)$, for finite and crown-free $T$}
\label{S:CFq-id}

Let $(\theta)$ be the following lattice-theoretical quasi-identity:
\begin{multline*}
\Bigl[a\leq(a'\vee b)\wedge(a'\vee c)\ \&\
b\leq(b'\vee a)\wedge(b'\vee c)\ \&\
c\leq(c'\vee a)\wedge(c'\vee b)\ \&\\
\&\ (a'\wedge a)\vee(b'\wedge b)\vee(c'\wedge c)\vee
(a\wedge b)\vee(a\wedge c)\vee(b\wedge c)\leq a'\wedge b'\wedge c'
\Bigr]\\
\Longrightarrow a\leq a'.
\end{multline*}
It is inspired by Example~\ref{Ex:SINonAt}
(see Corollary~\ref{C:CFq-id}).
The main result of Section~\ref{S:CFq-id} is the following.

\begin{theorem}\label{T:CFq-id}
Let $(T,\utr)$ be a finite crown-free poset. Then $\Co(T)$ satisfies
$(\theta)$.
\end{theorem}

Let us begin with an arbitrary (not necessarily finite, not
necessarily crown-free) poset $(T,\utr)$ and convex subsets $A$, $B$,
$C$, $A'$, $B'$, $C'$ of $T$ that satisfy the premise of $(\theta)$,
that is,
\begin{gather*}
 \begin{align*}
 A&\subseteq A'\vee B;& A&\subseteq A'\vee C;\\
 B&\subseteq B'\vee A;& B&\subseteq B'\vee C;\\
 C&\subseteq C'\vee A;& C&\subseteq C'\vee B;
 \end{align*}\\
 \begin{align*}
 A\cap A'&\subseteq B'\cap C';& B\cap B'&\subseteq A'\cap C';
 & C\cap C'&\subseteq A'\cap B';\\
 A\cap B&\subseteq A'\cap B';&
 B\cap C&\subseteq B'\cap C';&
 A\cap C&\subseteq A'\cap C'.
 \end{align*}
\end{gather*}
We shall put $\hA=A\setminus A'$, $\hB=B\setminus B'$, and
$\hC=C\setminus C'$. Observe that
 \begin{gather*}
 \hA\cap(B\cup C)=\hB\cap(A\cup C)=\hC\cap(A\cup B)=\es,\\
 \hA\cap\hB=\hA\cap\hC=\hB\cap\hC=\es.
 \end{gather*}
We shall later perform a construction whose key argument is provided by
the following lemma.

\begin{lemma}\label{L:OneLine}
Let $a\in\hA$ and let $a'\in A'$ with $a\tr a'$. Then there exists
$\seq{b,b'}\in\hB\times B'$ such that $b'\tr b\tr a$.
\end{lemma}

\begin{proof}
Observe first that $a\in A\subseteq A'\vee B$. Since
$a\notin A'\cup B$, there exists $(\ol{a}',b)\in A'\times B$ such
that either $\ol{a}'\tr a\tr b$ or $b\tr a\tr\ol{a}'$. In the first
case, $\ol{a}'\tr a\tr a'$, thus, by the convexity of $A'$,
$a\in A'$, a contradiction; whence
$b\tr a$. If $b\in B'$, then $b\in B\cap B'\subseteq A'$, but
$b\tr a\tr a'$, thus $a\in A'$, a contradiction; whence
$b\in\hB$. If there exists $x\in A$ with $x\utr b$, then, since
$b\tr a$, we obtain that $b\in A\cap B\subseteq A'$, a
contradiction again. But $b\in B\subseteq A\vee B'$ and $b\notin B'$,
thus there exists $b'\in B'$ such that $b'\tr b$.
\end{proof}

In particular, we observe the following corollary.

\begin{corollary}\label{C:OneLine}
The sets $\hA$, $\hB$, and $\hC$ are either simultaneously empty or
simultaneously nonempty.
\end{corollary}

\begin{proof}
If $\hA$ is nonempty, we pick $a\in\hA$. So $a\in A'\vee B$ while
$a\notin A'\cup B$, thus there is $(a',b)\in A'\times B$ such that
either $b\tr a\tr a'$ or $a'\tr a\tr b$. In the first case, we
apply Lemma~\ref{L:OneLine} to deduce that $\hB\neq\es$. In the second
case, we apply the dual of Lemma~\ref{L:OneLine} to reach the same
conclusion.
\end{proof}

Now we suppose that $\hA$ is nonempty, and we pick $a_0\in\hA$. As in
the proof of Corollary~\ref{C:OneLine}, there exists $a'_0\in A'$ such
that either $a_0\tr a'_0$ or $a'_0\tr a_0$; by replacing $\utr$ with
its dual if needed, we may assume without loss of generality that
$a_0\tr a'_0$.

By Lemma~\ref{L:OneLine}, there are $\seq{b_0,b'_0}\in\hB\times B'$
and $\seq{c_1,c'_1}\in\hC\times C'$ such that $b'_0\tr b_0\tr a_0$ and
$c'_1\tr c_1\tr a_0$. By applying the dual of Lemma~\ref{L:OneLine}
to $c'_1\tr c_1$, we obtain $\seq{b_1,b'_1}\in\hB\times B'$ such that
$c_1\tr b_1\tr b'_1$. By applying Lemma~\ref{L:OneLine} to
$b_1\tr b'_1$, we obtain $\seq{a_2,a'_2}\in\hA\times A'$ such that
$a'_2\tr a_2\tr b_1$. By applying in the same fashion
Lemma~\ref{L:OneLine} and its dual, we obtain
$\seq{c_2,c'_2}\in\hC\times C'$, $\seq{b_3,b'_3}\in\hB\times B'$, and
$\seq{a_3,a'_3}\in\hA\times A'$ such that $a_2\tr c_2\tr c'_2$,
$b'_3\tr b_3\tr c_2$, and $b_3\tr a_3\tr a'_3$.

Now we observe that $b'_0\tr b_0\tr a_0\tr a'_0$ and
$b'_3\tr b_3\tr a_3\tr a'_3$, that is, we can start the process
again. Arguing by induction, we obtain elements
$\seq{a_i,a'_i}\in\hA\times A'$ for $i\not\equiv 1\pmod{3}$, elements
$\seq{b_i,b'_i}\in\hB\times B'$ for $i\not\equiv 2\pmod{3}$, and
elements $\seq{c_i,c'_i}\in\hC\times C'$ for $i\not\equiv 0\pmod{3}$
such that the following relations hold, for all $i<\omega$:
 \begin{gather}
 b'_{3i}\tr b_{3i}\tr a_{3i}\tr a'_{3i};\label{Eq:b'a'}\\
 c'_{3i+1}\tr c_{3i+1}\tr b_{3i+1}\tr b'_{3i+1};\label{Eq:c'b'}\\
 a'_{3i+2}\tr a_{3i+2}\tr c_{3i+2}\tr c'_{3i+2}.\label{Eq:a'c'}
 \end{gather}
This can be illustrated by Figure~5.

\begin{figure}[htb]
\includegraphics{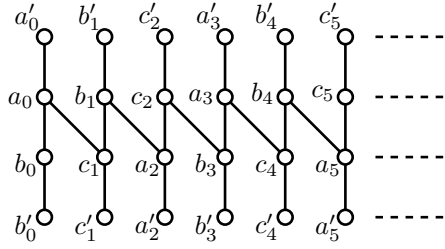}
\caption{A pattern in $T$.}
\end{figure}

Now we define subsets of $T$ as follows:
 \begin{align*}
 \Omega^+&=\setm{a_{3i}}{i<\omega}\cup\setm{b_{3i+1}}{i<\omega}
 \cup\setm{c_{3i+2}}{i<\omega};\\
 \Omega^-&=\setm{a_{3i+2}}{i<\omega}\cup\setm{b_{3i}}{i<\omega}
 \cup\setm{c_{3i+1}}{i<\omega};\\
 \Omega&=\Omega^+\cup\Omega^-.
 \end{align*}
Since $\hA$, $\hB$, and $\hC$ are mutually disjoint and their union
contains $\Omega$, we can define a map $\chi\colon\Omega\to 3$ by the
rule
 \[
 \chi(x)=\begin{cases}
 0&(x\in\hA),\\
 1&(x\in\hB),\\
 2&(x\in\hC),
 \end{cases}
 \qquad\text{for all }x\in\Omega.
 \]

\begin{lemma}\label{L:samechi}
For all $\seq{x,y}\in\Omega^-\times\Omega^+$, $\chi(x)=\chi(y)$
implies that $x\nutr y$. In particular, $\Omega^-\cap\Omega^+=\es$.
\end{lemma}

\begin{proof}
We need to prove that for all natural numbers $i$ and $j$, the
following inequalities hold:
\begin{itemize}
\item $a_{3i+2}\nutr a_{3j}$. Otherwise, by
\eqref{Eq:b'a'} and \eqref{Eq:a'c'}, $a'_{3i+2}\tr a_{3i+2}\tr a'_{3j}$,
thus $a_{3i+2}\in A'$, a contradiction.

\item $b_{3i}\nutr b_{3j+1}$. Otherwise, by
\eqref{Eq:b'a'} and \eqref{Eq:c'b'}, $b'_{3i}\tr b_{3i}\tr b'_{3j+1}$,
thus $b_{3i}\in B'$, a contradiction.

\item $c_{3i+1}\nutr c_{3j+2}$. Otherwise, by
\eqref{Eq:c'b'} and \eqref{Eq:a'c'},
$c'_{3i+1}\tr c_{3i+1}\tr c'_{3j+2}$, thus $c_{3i+1}\in C'$, a
contradiction.
\end{itemize}
This concludes the proof.
\end{proof}

For an integer $m\geq2$, we define a \emph{$m$-pre-crown} to be
a finite sequence $\famm{\seq{x_i,y_i}}{i\in\Zm}$ of elements of
$\Omega^-\times\Omega^+$ such that the following conditions hold, for
all $i\in\Zm$:
\begin{itemize}
\item[(C1)] $x_i,x_{i+1}\tr y_i$;

\item[(C2)] $\chi(x_i)\neq\chi(x_{i+1})$ and
$\chi(y_i)\neq\chi(y_{i+1})$ if $i\neq m-1$.
\end{itemize}

If $m=2$, then, by (C1), $x_0,x_1\tr y_0,y_1$. Furthermore, by (C2),
$\chi(x_0)\neq\chi(x_1)$, thus it follows from $x_0,x_1\tr y_0$ and
Lemma~\ref{L:samechi} that $\chi(y_0)$ is the unique element of
$3\setminus\set{\chi(x_0),\chi(x_1)}$. The same holds for $\chi(y_1)$,
whence $\chi(y_0)=\chi(y_1)$, which contradicts (C2). Therefore,
if there exists a $m$-pre-crown, then $m\geq3$.

We can now prove the main lemma of this section.

\begin{lemma}\label{L:cf2pcf}
Suppose that $T$ is crown-free. Then there are no pre-crowns in $T$.
\end{lemma}

\begin{proof}
Otherwise, let $m$ be the least positive integer such that there
exists a $m$-pre-crown, and let $\fc=\famm{\seq{x_i,y_i}}{i\in\Zm}$
be such a pre-crown. As observed before, $m\geq3$. By
assumption on $T$, in order to get a contradiction, it suffices to
prove that $\fc$ is a crown of $T$. By (C1) and
Lemma~\ref{L:subcrown}, it suffices to prove that for all $i$,
$j\in\Zm$ such that $i\notin\set{j,j+1}$, the inequality
$x_i\utr y_j$ does not hold. Suppose otherwise; by
Lemma~\ref{L:samechi}, $x_i\tr y_j$. Two cases can occur.

\begin{itemize}
\item[\textbf{Case 1.}] $i<j$. Then the finite sequence
 \[
 \seq{\seq{x_i,y_i},\seq{x_{i+1},y_{i+1}},\ldots,\seq{x_j,y_j}}
 \]
is a $(j-i+1)$-pre-crown (see Figure~6(i)),
with $1\leq j-i\leq m-1$. By the minimality
assumption on $m$, this cannot happen unless $i=0$ and $j=m-1$, in
which case $i=j+1$ (modulo $m$ as usual), a contradiction.

\begin{figure}[htb]
\includegraphics{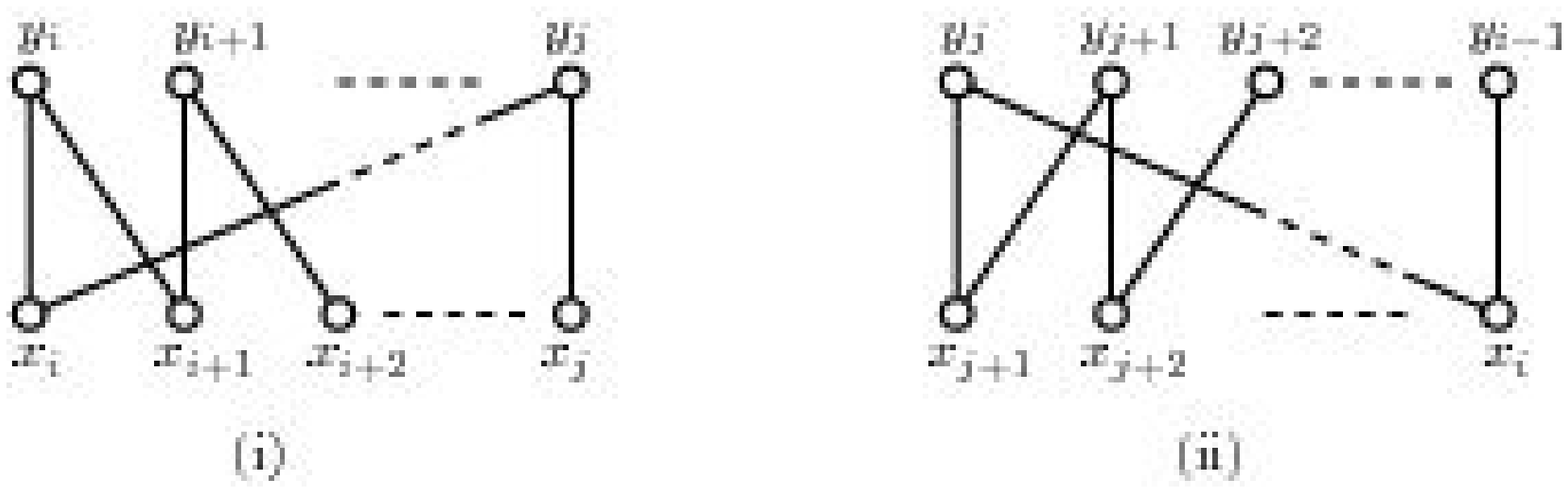}
\caption{Shorter pre-crowns.}
\end{figure}

\item[\textbf{Case 2.}] $j<i$. Then the finite sequence
 \[
 \seq{\seq{x_i,y_{i-1}},\ldots,\seq{x_{j+2},y_{j+1}},
 \seq{x_{j+1},y_j}}
 \]
is a $(i-j)$-pre-crown (see Figure~6(ii)),
with $2\leq i-j<m$, which contradicts again the minimality of $m$.
\end{itemize}

Hence $\fc$ is a $m$-crown of $T$, a contradiction.
\end{proof}

Now we have all the necessary tools to conclude the proof of
Theorem~\ref{T:CFq-id}.

\begin{proof}[Proof of Theorem~\textup{\ref{T:CFq-id}}]
Suppose that $T$ is finite and crown-free. There are $i<j$ such that
$b_{3i}=b_{3j}$. Then the finite sequence
 \[
 \seq{\seq{b_{3i},a_{3i}},\seq{c_{3i+1},b_{3i+1}},
 \ldots,\seq{a_{3j-1},c_{3j-1}}}
 \]
is a $(3j-3i)$-pre-crown in $T$ (see Figure~7), a contradiction.

\begin{figure}[htb]
\includegraphics{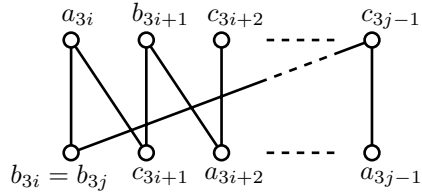}
\caption{A pre-crown in $T$.}
\end{figure}

Hence we have proved that $\hA=\es$, that is, $A\subseteq A'$.
Therefore, $\Co(T)$ satisfies~$(\theta)$.
\end{proof}

\begin{corollary}\label{C:CFq-id}
Let $Q$ be the finite poset and $\LL$ the finite lattice
of Example~\textup{\ref{Ex:SINonAt}}. Then, although $\LL$ embeds into
$\Co(Q)$, there is no finite, tree-like poset $R$ such that~$\LL$ embeds
into $\Co(R)$.
\end{corollary}

\begin{proof}
It follows from Proposition~\ref{P:TrCf} that $R$ is crown-free, thus,
by Theorem~\ref{T:CFq-id}, $\Co(R)$ satisfies $(\theta)$. On the other
hand, the lattice $\LL$ of Example~\ref{Ex:SINonAt} does not
satisfy~$(\theta)$ (consider the atoms $A$, $B$, $C$, $A'$, $B'$,
$C'$ of $\LL$), therefore it cannot be embedded into $\Co(R)$.
\end{proof}

On the other hand, it follows from Theorem~\ref{T:Main2}(iii) that if
a finite lattice $L$ without $\DD$-cycle embeds into some $\Co(P)$, then
it embeds into $\Co(R)$ for some finite, tree-like poset $R$. In the
presence of $\DD$-cycles anything can happen, for example, take
$L=\Co(\mathbf{4})$, the lattice of all order-convex subsets of a
four-element chain; it embeds into $\Co(\mathbf{4})$ for the
finite, tree-like poset $\mathbf{4}$, however it has $\DD$-cycles.

\section{Finite generation and word problem in $\SUB$}
\label{S:WordPb}

For a lattice term $\bs(\vx_1,\dots,\vx_n)$, a poset
$P$, and convex subsets $X_1$, \dots, $X_n$ of $P$, we denote by
$\bs^P(X_1,\dots,X_n)$ the evaluation of the term
$\bs(\vx_1,\dots,\vx_n)$ at $\seq{X_1,\dots,X_n}$ in the
lattice $\Co(P)$.

The present section rests on the following lemma. Its proof is an
easy induction argument on the length of $\bs$, that we leave to the
reader.

\begin{lemma}\label{L:sP(x)}
Let $n$ be a positive integer, let $\bs(\vx_1,\dots,\vx_n)$ be a
lattice term, and let $X_1$, \dots, $X_n$ be convex
subsets of a poset $P$. Then $\bs^P(X_1,\dots,X_n)$ is the directed
union of all subsets of the form $\bs^Q(X_1\cap Q,\dots,X_n\cap Q)$,
for $Q\subseteq P$ finite.
\end{lemma}

As immediate corollaries, we get the following.

\begin{corollary}\label{C:IdenPos}
Let $P$ be a poset. Any lattice-theoretical identity valid in all
$\Co(Q)$, for $Q$ a finite subset of $P$, is also valid in $\Co(P)$.
\end{corollary}

\begin{corollary}\label{C:IdenPos2}
A lattice-theoretical identity is valid in $\SUB$ if{f} it holds in
$\Co(P)$ for every \emph{finite} poset $P$.
\end{corollary}

Consequently, the variety $\SUB$ is generated by its finite members.
By using the results of J.\,C.\,C. McKinsey \cite{McKin}, we obtain
the following consequence.

\begin{corollary}\label{C:EqnDec}
The word problem in the variety $\SUB$ is decidable.
\end{corollary}

This means that it is decidable whether a given lattice identity
$\bs(\vx_1,\dots,\vx_m)=\bt(\vx_1,\dots,\vx_m)$ holds in all lattices
of the form $\Co(P)$. A closer look at the proof of
Lemma~\ref{L:sP(x)} shows that it is sufficient to verify whether the
given identity holds in all $\Co(P)$ for $|P|\leq n$,
where $n$ is the supremum of the lengths of the terms $\bs$ and~$\bt$.

\section{Open problems}\label{S:Pbs}
We know that the class $\SUB$ is generated, as a variety,
by its finite members (see Corollary~\ref{C:IdenPos2}). We also know
that any finite lattice in $\SUB$ can be embedded into some finite
$\Co(P)$ (see Theorem~\ref{T:Main}). Nevertheless we do not know
whether the latter generate the whole \emph{quasi}variety.

\begin{problem}\label{Pb:FinMem}
Is the class $\SUB$ generated, as a quasivariety, by its finite
members?
\end{problem}

Equivalently, does there exist a lattice
quasi-identity that holds in all finite $\Co(P)$-s but not in all
$\Co(P)$-s?

\begin{problem}\label{Pb:UTDec}
Is the universal theory of all lattices of the form $\Co(P)$
decidable?
\end{problem}

A positive answer to Problem~\ref{Pb:FinMem} would yield a positive
answer to Problem~\ref{Pb:UTDec}.

\begin{problem}\label{Pb:Chains}
Is the class $\mathbf{C}$ of all lattices that can be embedded into a
product of the form $\prod_{i\in I}\Co(C_i)$, where the $C_i$ are
\emph{chains}, a variety?
\end{problem}

Problem~\ref{Pb:Chains} is solved by the authors in \cite{SeWe3}.

\begin{problem}\label{Pb:Functor}
Can the embedding problem of a lattice in $\SUB$
into some $\Co(P)$ be solved by a \emph{functor} (that, say, sends any
$L$ to some $\Co(P)$)? Can such a functor be idempotent?
\end{problem}

Our next problem has a more computational nature.

\begin{problem}\label{Pb:CardP}
For each positive integer $n$, denote by $\xi(n)$ the least positive
integer such that every finite lattice $L$ in $\SUB$ with $n$ \jirr\
elements embeds into some $\Co(P)$, where
$|P|\leq\xi(n)$. Compute $\xi(n)$, for all $n>0$. Does $\xi(n)=O(n)$
as~$n$ goes to infinity?
\end{problem}

For a sublattice $K$ of a finite lattice $L$, the inequality
$|\J(K)|\leq|\J(L)|$ holds, see~\cite[Lemma~2]{Ad1}.
In particular, if a finite lattice $L$ embeds into $\Co(P)$ for some
finite poset~$P$, then $|\J(L)|\leq|P|$. By combining this with the
result of Theorem~\ref{T:Main}, we obtain the inequalities
 \[
 n\leq\xi(n)\leq 2n^2-5n+4.
 \]

\section*{Acknowledgments}
This work was started during the second author's visit at the
Institute of Mathematics of the Siberian Branch of RAS, in Novosibirsk,
during the summer of 2001. The outstanding hospitality met there is
greatly appreciated.

The authors wish to thank warmly Kira Adaricheva and Mikhail Sheremet
for their interest and suggestions about the present work. This work
was completed while both authors were visiting the Charles University
during the fall of 2001. Excellent conditions provided by the
Department of Algebra are highly appreciated. Special thanks are due to
Ji\v{r}\'\i\ T\r{u}ma and V\'aclav Slav\'\i k.

\end{document}